\documentclass[11pt]{article}%
\usepackage{amssymb}
\usepackage[all]{xy}
\usepackage[latin1]{inputenc}
\usepackage{amsfonts}
\usepackage{amsmath}
\usepackage{amssymb}
\usepackage{url}
\usepackage{hyperref}
\usepackage{graphicx}%
\providecommand{\U}[1]{\protect\rule{.1in}{.1in}}
\numberwithin{equation}{section}
\providecommand{\U}[1]{\protect\rule{.1in}{.1in}}
\providecommand{\U}[1]{\protect\rule{.1in}{.1in}}
\newtheorem{theo}{Theorem}[section]
\newtheorem{prop}[theo]{Proposition}
\newtheorem{lem}[theo]{Lemma}
\newtheorem{cor}[theo]{Corollary}

\newtheorem{rem}[theo]{Remark}

\begin{document}

\title{A Lognormal Central Limit Theorem for Particle Approximations of Normalizing Constants}
\author{Jean B\'erard\thanks{Universit\'e Claude Bernard, Lyon 1, France.
\ \texttt{jean.berard@univ-lyon1.fr}} , Pierre Del Moral\thanks{Centre INRIA
Bordeaux Sud-Ouest Institut de Math\'ematiques, Universit\'e Bordeaux I, 33405
Talence cedex, France.\ \texttt{pierre.del-moral@inria.fr.}} , Arnaud
Doucet\thanks{Department of Statistics, Oxford
University.\ \texttt{doucet@stats.ox.ac.uk.}}}
\maketitle

\begin{abstract}
This paper deals with the numerical approximation of normalizing constants
produced by particle methods, in the general framework of Feynman-Kac
sequences of measures. It is well-known that the corresponding estimates
satisfy a central limit theorem for a fixed time horizon $n$ as the number of
particles $N$ goes to infinity. Here, we study the situation where both $n$
and $N$ go to infinity in such a way that $\lim_{n\rightarrow\infty}%
n/N=\alpha>0$. In this context, Pitt et al. \cite{pitt2012} recently
conjectured that a lognormal central limit theorem should hold. We formally
establish this result here, under general regularity assumptions on the model.
We also discuss special classes of models (time-homogeneous environment and
ergodic random environment) for which more explicit descriptions of the
limiting bias and variance can be obtained.


\emph{Mathematics Subject Classification} : \newline Primary: 65C35, 47D08,
60F05; Secondary: 65C05, 82C22, 60J70.

\end{abstract}

\section{Introduction}

\subsection{Feynman-Kac measures and their particle approximations}

Consider a Markov chain $(X_{n})_{n\geq0}$ on a measurable state space
$(E,\mathcal{E})$, whose transitions are prescribed by a sequence of Markov
kernels $\left(  M_{n}\right)  _{n\geq1}$, and a collection of positive
bounded and measurable functions $(G_{n})_{n\geq0}$ on $E$. We associate to
$\left(  M_{n}\right)  _{n\geq1}$ and $(G_{n})_{n\geq0}$ the sequence of
unnormalized Feynman-Kac measures $\left(  \gamma_{n}\right)  _{n\geq0}$ on
$E$, defined through their action on bounded (real-valued) measurable
functions by:
\begin{equation}
\gamma_{n}(f):=\mathbb{E}\left(  f(X_{n})~\prod_{0\leq p<n}G_{p}
(X_{p})\right)  . \label{FK-u}%
\end{equation}
The corresponding sequence of normalized (probability) Feynman-Kac measures
$\left(  \eta_{n}\right)  _{n\geq0}$ is defined by:
\begin{equation}
\eta_{n}(f):={\gamma_{n}(f)}/{\gamma_{n}(1)} . \label{FK-n}%
\end{equation}
It is easily checked that, for all $n \geq0$, the normalizing constant
$\gamma_{n}(1)$ satisfies
\begin{equation}
\gamma_{n}(1)=\mathbb{E}\left(  \prod_{0\leq p<n}G_{p}(X_{p})\right)
=\prod_{0\leq p<n}\eta_{p}(G_{p}). \label{multiplicative}%
\end{equation}
Here and throughout the paper, the notation $\mu(f)$, where $\mu$ is a
finite signed measure and $f$ is a bounded function defined on the same space, is
used to denote the Lebesgue integral of $f$ with respect to $\mu$, i.e.
$\mu(f) := \int f(x) d\mu(x)$. Given a bounded integral operator $K(x,dx^{\prime})$ from $E$ into itself, we denote by
$\mu K$ the measure resulting from the action of $K$ on $\mu$,
i.e.
\[
\mu K (dx^{\prime}) := \int\mu(dx) K(x, dx^{\prime}).
\]
For a bounded measurable function $f$ on $E$, we denote by $K(f)$ the (bounded
measurable) function resulting from the action of $K$ on $f$, i.e.
\[
K(f)(x) := \int f(x^{\prime}) K(x, dx^{\prime}).
\]

Feynman-Kac measures appear in numerous scientific fields including, among
others, signal processing, statistics and statistical physics; see
\cite{cappe2}, \cite{DMo04} and \cite{doucet2001} for many applications. For
example, in a non-linear filtering framework, the measure $\eta_{n}$
corresponds to the posterior distribution of the latent state of a dynamic
model at time $n$ given the observations collected from time $0$ to time
$n-1$, and $\gamma_{n}(1)$ corresponds to the likelihood of these very
observations. A generic Monte Carlo application has $(\eta_{n})_{n\geq0}$
corresponding to a sequence of tempered versions of a distribution $\eta$ that
we are interested in sampling from using suitable $\eta_{n}$-invariant Markov
kernels $M_{n}$, with $(\gamma_{n}(1))_{n\geq0}$ the resulting sequence of
normalizing constants \cite{DelmoralDoucetJasra2006}. Two applications are
discussed in more details in Sections~\ref{part-absorp}
and~\ref{filtering-sec}.

A key issue with Feynman-Kac measures is that they are analytically
intractable in most situations of interest. Over the past twenty years,
particle methods have emerged as the tool of choice to produce numerical
approximations of these measures and their associated normalizing constants.
We give a brief overview of these methods here, and refer to \cite{DMo04} for
a more thorough treatment.

We first observe that the sequence $(\eta_{n})_{n \geq0}$ admits the following
inductive representation: for all $n \geq1$, one has
\begin{equation}
\eta_{n} = \Phi_{n}(\eta_{n-1}). \label{rec-base}%
\end{equation}
Here, $\Phi_{n}$ is the non-linear transformation on probability measures
defined by
\[
\Phi_{n}(\mu):= \Psi_{G_{n-1}}(\mu)M_{n},
\]
where, given a bounded positive function $G$ and a probability measure $\mu$
on $E$, $\Psi_{G}$ denotes the Boltzmann-Gibbs transformation:
\begin{equation}
\Psi_{G}(\mu)(dx):=\frac{1}{\mu(G)}G(x)\mu(dx) . \label{def-Boltzmann-Gibbs}%
\end{equation}

One then looks for representations of $\Phi_{n}$ of the form:
\begin{equation}
\Phi_{n}(\mu) = \mu K_{n,\mu} , \label{phidef}%
\end{equation}
where $(K_{n,\mu})_{n , \mu}$ is a collection of Markov kernels defined for
every time-index $n \geq1$ and probability measure $\mu$ on $E$. The choice
for $K_{n,\mu}$ is far from being unique. One can obviously use $K_{n,\mu
}(x,dx^{\prime}):=\Phi_{n}(\mu)(dx^{\prime})$, but there are alternatives. For
example, if $G_{n-1}$ takes its values in the interval $]0,1]$, $\Psi
_{G_{n-1}}(\mu)$ can be expressed through a non-linear Markov transport
equation
\begin{equation}
\Psi_{G_{n-1}}(\mu)=\mu S_{G_{n-1},\mu} \label{nonlin-transport}%
\end{equation}
with the non-linear Markov transition kernel
\[
S_{G_{n-1},\mu}(x,dx^{\prime}):=G_{n-1}(x)~\delta_{x}(dx^{\prime})+\left(
1-G_{n-1}(x)\right)  ~\Psi_{G_{n-1}}(\mu)(dx^{\prime})
\]
so we can use
\begin{equation}
K_{n,\mu}:=S_{G_{n-1},\mu}M_{n}. \label{2nd-choice}%
\end{equation}

The non-linear Markov representation (\ref{phidef}) directly suggests a
mean-field type particle approximation scheme for $\left(  \eta_{n}\right)
_{n\geq0}$. For every $n\geq0$, we have an $N-$tuple of elements of $E$
denoted by $\xi_{n}^{(N)}=\left(  \xi_{n}^{(N,i)}\right)  _{1\leq i\leq N}$,
whose empirical measure $\eta_{n}^{N}:=\frac{1}{N}\sum_{j=1}^{N}\delta
_{\xi_{n}^{(N,j)}}$ provides a particle approximation of $\eta_{n}$. The
sequence $(\xi_{n}^{(N)})_{n\geq0}$ evolves as an $E^{N}$-valued Markov chain
whose initial distribution is given by $\mathbb{P}\left(  \xi_{0}^{(N)}\in
dx\right)  =\prod_{i=1}^{N}\eta_{0}\left(  dx_{i}\right)  $, while, for
$n\geq1$, the transition mechanism is specified by
\begin{equation}
\text{ }\mathbb{P}\left(  \xi_{n}^{(N)}\in dx~\left\vert ~\mathcal{F}%
_{n-1}^{N}\right.  \right)  =\prod_{i=1}^{N}~K_{n,\eta_{n-1}^{N}}(\xi
_{n-1}^{(N,i)},dx^{i}). \label{meanfield}%
\end{equation}
Here $\mathcal{F}_{n-1}^{N}$ is the sigma-field generated by the random
variables $(\xi_{p}^{(N)})_{0\leq p\leq n-1}$, and $dx:=dx^{1}\times
\ldots\times dx^{N}$ stands for an infinitesimal neighborhood of a point
$x=(x^{1},\ldots,x^{N})\in E^{N}$.

Using the identity (\ref{multiplicative}) we can easily obtain a particle
approximation $\gamma_{n}^{N}\left(  1\right)  $ of the normalizing constant
$\gamma_{n}\left(  1\right)  $ by replacing the measures $\left(  \eta
_{p}\right)  _{p=0}^{n-1}$ by their particle approximations $\left(  \eta
_{p}^{N}\right)  _{p=0}^{n-1}$ to get
\begin{equation}
\gamma_{n}^{N}(1):=\prod_{0\leq p<n}\eta_{p}^{N}\left(  G_{p}\right)
\label{particleapproxnormalizingconstant}%
\end{equation}
and we define its normalized version by
\begin{equation}
\overline{\gamma}_{n}^{N}(1)={\gamma_{n}^{N}(1)}/{\gamma_{n}(1)}=\prod_{0\leq
p<n}\eta_{p}^{N}(\overline{G}_{p})\quad\mbox{\rm with}\quad\overline{G}%
_{n}:={G_{n}}/{\eta_{n}(G_{n}).} \label{relativeestimate}%
\end{equation}

The main goal of this article is to establish a central limit theorem for
$\log$ $\overline{\gamma}_{n}^{N}(1)$ as $n\rightarrow\infty$ when the number
of particles $N$ is proportional to $n$. Such a result has been conjectured by
Pitt et al. \cite{pitt2012}, who provided compelling empirical evidence for
it. To our knowledge, the present work gives the first mathematical proof of a
result of this type.

\subsection{Statement of the main result\label{section:statementresults}}

To state our result, we need to introduce additional notations. We start with
the convention that $\Phi_{0}(\mu) := \eta_{0}$ for all $\mu$, $K_{0,\mu}(x,
\cdot):=\eta_{0}(\cdot)$ for all $x$, and $\mathcal{F}_{-1}^{N}=\{\emptyset
,\Omega\}$. For the sake of definiteness, we also let $\eta_{-1} := \eta_{0}$
and $\eta_{-1}^{N} := \eta_{0}$. These conventions make
\eqref{rec-base}-\eqref{phidef}-\eqref{meanfield} valid for $n=0$.

Then denote by $V_{n}^{N}$ the centered local error random fields defined, for
$n \geq0$, by
\begin{equation}
V_{n}^{N}:=\sqrt{N}\text{ }\left(  \eta_{n}^{N}-\Phi_{n}(\eta_{n-1}
^{N})\right)  , \label{VnNlocalsamplingerror}%
\end{equation}
so that one can write
\[
\eta_{n}^{N}=\Phi_{n}(\eta_{n-1}^{N})+\frac{1}{\sqrt{N}}V_{n}^{N}.
\]

To describe the corresponding covariance structure, let us introduce, for all
$n\geq0$, bounded functions $f_{1},f_{2}$, and probability measure $\mu$, the
notation
\[
\mbox{\rm Cov}_{n,\mu}(f_{1},f_{2}):=\mu\left[  K_{n,\mu}(f_{1}f_{2}%
)-K_{n,\mu}(f_{1})K_{n,\mu}(f_{2})\right]  .
\]
We then have the following explicit expression for conditional covariances:
\begin{equation}
\mathbb{E}\left(  V_{n}^{N}(f_{1})V_{n}^{N}(f_{2})\left\vert \mathcal{F}%
_{n-1}^{N}\right.  \right)  =\mbox{\rm Cov}_{n,\eta_{n-1}^{N}}(f_{1},f_{2}).
\label{def-Cov-eta}%
\end{equation}
It is proved in~\cite[chapter 9]{DMo04} that, under weak regularity
assumptions, $(V_{n}^{N})_{n\geq0}$ converges in law, as $N$ tends to
infinity, to a sequence of $n$ independent, Gaussian and centered random
fields $(V_{n})_{n\geq0}$ with a covariance given by
\begin{equation}
C_{V_{n}}(f_{1},f_{2}):=\mathbb{E}(V_{n}(f_{1})V_{n}(f_{2}%
))=\mbox{\rm Cov}_{n,\eta_{n-1}}(f_{1},f_{2}). \label{corr1}%
\end{equation}
Note that, with the special choice $K_{n,\mu}(x,\cdot):=\Phi_{n}(\mu)$,
\eqref{corr1} reduces to
\begin{equation}
C_{V_{n}}(f_{1},f_{2})=\eta_{n}(f_{1}f_{2})-\eta_{n}(f_{1})\eta_{n}(f_{2}).
\label{corrspecialcase}%
\end{equation}

Let us now introduce the family of operators $(Q_{p,n})_{0 \leq p \leq n}$
acting on the space of bounded measurable functions, defined by
\begin{equation}
\quad Q_{p,n}(f)(x):=\mathbb{E} \left(  \left.  f(X_{n})~\prod_{p\leq
q<n}G_{q}(X_{q}) \right|  X_{p}=x \right)  . \label{Qpn}%
\end{equation}
It is easily checked that $(Q_{p,n})_{0 \leq p \leq n}$ forms a semigroup for
which $\gamma_{n}=\gamma_{p}Q_{p,n}$.

We also define
\begin{equation}
\overline{Q}_{p,n}(f):=\frac{{Q_{p,n}(f)}}{{\eta_{p}Q_{p,n}(1)}}.
\label{eq:normalizedQpn}%
\end{equation}
Finally, we define the Markov kernel $P_{p,n}$ through its action on bounded
measurable functions:
\begin{equation}
P_{p,n}(f):=Q_{p,n}(f)/Q_{p,n}(1). \label{Ppn}%
\end{equation}

It is well-known in the literature that (see for example \cite[chapter
9]{DMo04}), for fixed $n$, as $N \to+\infty$, the following convergence in
distribution holds under weak regularity assumptions:
\begin{equation}
\sqrt{N} \left(  \overline{\gamma}_{n}^{N}(1) -1 \right)
\xrightarrow[N \to +\infty]{d} \sum_{0\leq p<n}V_{p}(\overline{Q}_{p,n}(1)).
\label{CLT-classique}%
\end{equation}

Here, we are here interested in the fluctuations of $\overline{\gamma}_{n}%
^{N}(1)$ as both $n,N\rightarrow\infty$ with $N$ proportional to $n$. It turns
out that, in such a regime, the observed behavior is different from that
described by \eqref{CLT-classique}. Indeed, the magnitude of the fluctuations
of $\overline{\gamma}_{n}^{N}(1)$ around $1$ does not vanish as $n,N$ go to
infinity, and they are described in the limit by a log-normal instead of a
normal distribution.

Our result is obtained under specific assumptions that we now list. First, the
potential functions are assumed to satisfy
\begin{equation}
\label{H0}g_{n}:={\sup G_{n}}/{\inf G_{n}}<+\infty\quad\mbox{\rm and}\quad
g:=\sup_{n\geq0}g_{n}<+\infty.
\end{equation}
Moreover, we assume that the Dobrushin coefficient of $P_{p,n}$, denoted
$\beta(P_{p,n})$, satisfies
\begin{equation}
\beta(P_{p,n})\leq a~e^{-\lambda(n-p)} \label{H}%
\end{equation}
for some finite constant $a<+\infty$ and some positive $\lambda>0$. Finally,
we assume that the kernels $K_{n, \mu}$ satisfy an inequality of the following
form:
\begin{equation}
\Vert\left[  K_{n,\mu_{1}}-K_{n,\mu_{2}}\right]  (f)\Vert\leq\kappa~|(\mu_{1}
-\mu_{2})(T_{n}(f, \mu_{2}))|, \label{H2}%
\end{equation}
for any two probability measures $\mu_{1}, \mu_{2}$ on $E$, and any measurable
map $f$ with oscillation $\mbox{\rm osc}(f):=\mbox{\rm Sup}_{x,y} |f(x)-f(y)|
\leq1$, where $\kappa$ is a finite constant, and $T_{n}(f, \mu_{2})$ is a
measurable map with oscillation $\leq1$ that may depend on $n,f, \mu_{2}$.

In the rest of the paper, unless otherwise stated, we assume that
(\ref{H0})-(\ref{H})-(\ref{H2}) hold.

Several sufficient conditions on the Markov kernels $M_{n}$ under which
(\ref{H}) holds are discussed in~\cite[Section 4.3]{DMo04}, as well as in
Section 3.4 in~\cite{dhw-2012}. Conditions under which (\ref{H2}) is satisfied
are given in Section \ref{cov-regularity}.

We are now in position to state the main result of the paper.

\begin{theo}
\label{lognormal} Assume (\ref{H0})-(\ref{H})-(\ref{H2}), and let $v_{n}$ be
defined as
\[
v_{n}:=\sum_{0\leq p<n}\mathbb{E}\left(  V_{p}(\overline{Q}_{p,n}%
(1))^{2}\right)  =\sum_{0\leq q<n}\mbox{\rm Cov}_{q,\eta_{q-1}}(\overline
{Q}_{q,n}(1),\overline{Q}_{q,n}(1)).
\]
Assume that $N$ depends on $n$ in such a way that
\[
\lim_{n\rightarrow+\infty}\frac{n}{N}=\alpha\in]0,+\infty[,
\]
and that
\begin{equation}
\lim_{n\rightarrow+\infty}\frac{v_{n}}{n} = \sigma^{2} \in]0,+\infty[.
\label{assumption-borneinf}%
\end{equation}
One then has the following convergence in distribution:
\begin{equation}
\log{\overline{\gamma}_{n}^{N}(1)} \xrightarrow[n \to +\infty]{d}
\mathcal{N}\left(  -\frac{1}{2}\alpha\sigma^{2},\alpha\sigma^{2}\right)  ,
\label{eq:asymptoticlaw}%
\end{equation}
where $\mathcal{N}\left(  u,v \right)  $ denotes the normal distribution of
mean $u$ and variance $v$.
\end{theo}

\begin{rem}
It follows\ from the continuous mapping theorem that ${\overline{\gamma}%
_{n}^{N}(1)}$ asymptotically exhibits a log-normal distribution. The
relationship between the asymptotic bias and variance in
(\ref{eq:asymptoticlaw}) should not be a surprise since $\mathbb{E}%
(\overline{\gamma}_{n}^{N}(1))=1$ for any $n,N$ \cite[Proposition
7.4.1.]{DMo04}.
\end{rem}

\begin{rem}
We believe that Theorem \ref{lognormal} may be established under the weaker
stability assumptions developed in \cite{douc2012} and \cite{nickstab}, at the
price of a significantly increased technical complexity.
\end{rem}

\begin{rem}
Under assumption (\ref{H0}), it is easily seen that one always has $\sup_{n}
\frac{v_{n}}{n} < +\infty$. If, in addition to (\ref{H0})-(\ref{H}%
)-(\ref{H2}), one assumes that $\liminf_{n \to+\infty} \frac{v_{n}}{n}>0$
instead of the stronger assumption \eqref{assumption-borneinf}, the proof of
Theorem \ref{lognormal} still leads to a lognormal limit theorem of the
following form:
\[
\frac{1}{\sqrt{\alpha{{\mbox{\Large $\frac{v_n}{n}$}}}}} \left(
\log{\overline{\gamma}^{N}_{n}(1)} + \frac{\alpha}{2}\frac{v_{n}}{n} \right)
\xrightarrow[n \to +\infty]{d} \mathcal{N }(0,1).
\]

\end{rem}

This theoretical result was used in \cite{pitt2012} to optimize the asymptotic
variance of Metropolis-Hastings estimates, for a given computational budget,
using proposal distributions based on particle methods. Another
straightforward application is to the bias-correction of log-Bayes factors
estimates in large datasets. Yet another potential application in the spirit
of \cite{nick} is that $\sigma^{2}$ provides a criterion which could be used
to select between various interacting particle schemes.

\subsection{Some illustrations}

Here, we discuss two concrete situations where Theorem \ref{lognormal} can be
used, and where the variance expression (\ref{assumption-borneinf}) can be
made more explicit.

\subsubsection{Particle absorption models}

\label{part-absorp}Consider a particle in an absorbing random medium, whose
successive states $\left(  X_{n}\right)  _{n\geq0}$ evolve according to a
Markov kernel $M$. At time $n$, the particle is absorbed with probability
$1-G\left(  X_{n}\right)  $, where $G$ is a $\left[  0,1\right)  $-valued
potential function. Letting $G_{n}:=G$ for all $n \geq0$, and $M_{n} := M$ for
all $n \geq1$, the connection with the Feynman-Kac formalism is the following:
denoting by $T$ the absorption time of the particle, we have that $\gamma
_{n}(1)=\mathbb{P}\left(  T\geq n\right)  $, and $\eta_{n}%
=\mbox{\rm Law}\left(  X_{n}~|~T\geq n\right)  $. In this situation, the
multiplicative formula (\ref{multiplicative}) takes the form
\[
\mathbb{P}\left(  T\geq n\right)  =\prod_{0\leq m<n}\mathbb{P}\left(  T\geq
m+1 ~|~T\geq m\right)  ,
\]
where
\[
\mathbb{P}\left(  T\geq m+1~|~T\geq m\right)  =\int G(x)~\mathbb{P}\left(
X_{m}\in dx~|~T\geq m\right)  =\eta_{m}(G).
\]
In the present context, we have a map $\Phi$ such that $\Phi_{n} = \Phi$ for
all $n \geq1$, and conditions (\ref{H0})-(\ref{H}) ensure that $\Phi$ has a
unique fixed point measure $\eta_{\infty}$ such that
\[
\mbox{\rm Law}\left(  X_{n}\left\vert T\geq n\right.  \right)  \longrightarrow
_{_{n\to\infty}}\eta_{\infty}=\Phi(\eta_{\infty}).
\]
Moreover, we have that
\[
\overline{Q}_{0,n}(1)(x)={\mathbb{P}\left(  T\geq n\left\vert X_{0}=x\right.
\right)  }/{\mathbb{P}\left(  T\geq n\right)  }\longrightarrow_{_{_{n\to
\infty}}}h(x).
\]
Setting $\overline{Q}=Q/\eta_{\infty}Q(1)$, we find that the function $h$
satisfies the spectral equations
\[
\overline{Q}(h)=h\Leftrightarrow Q(h)=\lambda h, \mbox{\rm with } \lambda
=\eta_{\infty}(G).
\]

The measure $\eta_{\infty}$ is the so-called quasi-invariant or Yaglom
measure. Under some additional conditions, the parameter $\lambda$ coincides
with the largest eigenvalue of the integral operator $Q$, and $h$ is the
corresponding eigenfunction. In statistical physics, $Q$ comes from a
discrete-time approximation of a Schr\"odinger operator, and $h$ is called the
ground state function. For a more thorough discussion, we refer the reader to
Chapters 2 and 3 in~\cite{DMo04} and Chapter 7 in~\cite{dp-2013}.

In this scenario, the limiting variance $\sigma^{2}$ appearing in
(\ref{eq:asymptoticlaw}) is given by
\begin{equation}
\sigma^{2}=\mbox{\rm Cov}_{1,\eta_{\infty}}(h,h).
\label{eq:variancehomogeneous}%
\end{equation}
In particular, if the Markov kernels used in the particle approximation scheme
are given by $K_{\eta}(x,\mbox{\LARGE .})=\Phi(\eta)$, then using
(\ref{corrspecialcase}) we find that $\sigma^{2}=\eta_{\infty}\left(
[h-1]^{2}\right)  $. The detailed statement and proof of these results are
provided in Section~\ref{time-homog}.

\subsubsection{Non-linear filtering}

\label{filtering-sec}Let $(X_{n},Y_{n})_{n\geq0}$ be a Markov chain on some
product state space $E_{1}\times E_{2}$ whose transition mechanism takes the
form
\[
\mathbb{P}\left(  (X_{n},Y_{n})\in d(x,y)~|~(X_{n-1},Y_{n-1})\right)
=M_{n}(X_{n-1},dx)~g_{n}(y,x)~\nu_{n}(dy),
\]
where $\left(  \nu_{n}\right)  _{n\geq0}$ is a sequence of positive measures
on $E_{2}$, $\left(  M_{n}\right)  _{n\geq0\text{ }}$ is a sequence of Markov
kernels from $E_{1}$ into itself, and $\left(  g_{n}\right)  _{n\geq0}$ is a
sequence of density functions on $E_{2}\times E_{1}$. The aim of non-linear
filtering is to infer the unobserved process $\left(  X_{n}\right)  _{n\geq0}$
given a realization of the observation sequence $Y=y$. It is easy to check
that
\[
\eta_{n}=\mbox{\rm Law}\left(  X_{n}~|~Y_{m}=y_{m}~,~\forall0\leq m<n\right)
,
\]
using $G_{n}:=g_{n}(y_{n},\mbox{\LARGE .})$ in (\ref{FK-u}). Furthermore, the
density denoted $p_{n}(y_{0},\ldots,y_{n})$ of the random sequence of
observations $(Y_{0},\ldots,Y_{n})$ w.r.t. to the product measure
$\otimes_{0\leq p\leq n}\nu_{p}$ evaluated at the observation sequence, that
is the marginal likelihood, is equal to the normalizing constant $\gamma
_{n+1}(1)$. In this context, the multiplicative formula (\ref{multiplicative})
takes the following form
\[
p_{n}(y_{0},\ldots,y_{n})=\prod_{0\leq m\leq n}q_{m}(y_{m}~|~y_{l},~0\leq
l<m)
\]
with
\[
q_{m}(y_{m}~|~y_{l},~0\leq l<m)=\int g_{m}(y_{m},x)~\mathbb{P}\left(  X_{m}\in
dx~|~Y_{l}=y_{l},~0\leq l<m\right)  =\eta_{m}(G_{m}).
\]
For time-homogeneous models $(g_{m},M_{m})=(g,M)$ associated to an ergodic
process $Y$ satisfying a random environment version of Assumption (\ref{H}),
the ergodic theorem implies that the normalized log-likelihood function
converges to the entropy of the observation sequence
\begin{align*}
\frac{1}{n+1}\log{p_{n}(Y_{0},\ldots,Y_{n})}  &  =\frac{1}{n+1}\sum_{0\leq
m\leq n}\log{q}_{m}{(Y_{m}~|~Y_{l},~0\leq l<m)}\\
&  \longrightarrow_{n\to\infty}~~\mathbb{E}\left(  \log q(Y_{0}~|~Y_{m}%
,~m<0)\right)  ,
\end{align*}
where $q(Y_{0}~|~Y_{m},~m<0)$ is the conditional density of the random
variable $Y_{0}$ w.r.t. the infinite past. In
Section~\ref{FK-random-environment}, we shall prove the existence of a
limiting measure $\eta_{\infty}^{Y}$, and function $h^{Y}$ such that
\[
q(Y_{0}~|~Y_{m},~m<0)=\eta_{\infty}^{Y}(g(Y_{0},\mbox{\LARGE .}))
\]
and
\[
\overline{Q}_{0,n+1}^{Y}(1)(x):=\frac{q_{0,n}((Y_{0},\ldots,Y_{n})~|~x)}%
{\int\eta_{\infty}^{Y}(dx)~q_{0,n}((Y_{0},\ldots,Y_{n})~|~x)}\longrightarrow
_{n\to\infty}h^{Y}(x)
\]
where $q_{0,n}((Y_{0},\ldots,Y_{n})|x)$ stands for the conditional density of
$(Y_{0},\ldots,Y_{n})$ given $X_{0}=x$. Similar type results have been
recently established in \cite{nick} using slightly more restrictive
assumptions. In this situation, the limiting variance $\sigma^{2}$ appearing
in (\ref{eq:asymptoticlaw}) satisfies
\begin{equation}
\sigma^{2}=\mathbb{E}\left(  \mbox{\rm Cov}^{\theta^{-1}(Y)}_{1,\eta_{\infty
}^{\theta^{-1}(Y)}}(h^{Y},h^{Y}))\right)  , \label{eq:variancerandom}%
\end{equation}
where $\theta$ denotes the shift operator, and, if the Markov kernels used by
the particle approximation scheme are given by $K_{n,\eta}%
(x,\mbox{\LARGE .})=\Phi_{n}(\eta)$ associated to the potential $G_{n}%
:=g_{n}(Y_{n},\mbox{\LARGE .})$, then using (\ref{corrspecialcase}) we obtain%
\[
\sigma^{2}=\mathbb{E}\left(  \eta_{\infty}^{\theta^{-1}(Y)}\left(  \left[
h^{Y}-1\right]  ^{2}\right)  \right)  .
\]

The detailed statement and proof of these results are provided in
Section~\ref{FK-random-environment}.

\subsection{Notations and conventions}

We denote, respectively, by $\mathcal{M}(E)$, $\mathcal{P}(E)$ and
$\mathcal{B}_{b}(E)$, the set of all finite signed measures on space
$(E,\mathcal{E})$ equipped with total variation norm $\Vert
\mbox{\LARGE .}\Vert_{{\tiny \mathrm{tv}}}$, the subset of all probability
measures, and the Banach space of all bounded and measurable functions $f$
equipped with the uniform norm $\Vert f\Vert=\hbox{Sup}_{x\in E}|f(x)|$. We
also denote by $\mbox{\rm Osc}(E)$, the set of $\mathcal{E}$-measurable
functions $f$ with oscillations $\mbox{\rm osc}(f):=\mbox{\rm Sup}_{x,y}%
|f(x)-f(y)|\leq1$. We also denote by $\Vert X\Vert_{m}=\mathbb{E}%
(|X|^{m})^{1/m}$, the $\mathbb{L}_{m}$-norm of the random variable $X$, where
$m \geq1$.

In the sequel, the generic notation $c$ is used to denote a constant that
depends only on the model. To alleviate notations, we do not use distinct
indices (e.g. $c_{1}, c_{2}, \ldots$) each time such a constant appears, and
keep using the notation $c$ even though the corresponding constant may vary
from one statement to the other. Still, to avoid confusion, we sometimes make
a distinction between such constants by using $c,c^{^{\prime}},c^{^{\prime
\prime}}$ inside an argument. When the constant also depends on additional
parameters $p_{1},\ldots, p_{\ell}$, this is explicitly stated in the notation
by writing $c(p_{1},\ldots, p_{\ell})$.

\subsection{Organization of the paper}

The rest of the paper is organized as follows. In Section \ref{cov-regularity}%
, we establish basic regularity properties of the $\mathrm{Cov}$ operator.
Section \ref{sec-quantitative} is devoted to the long-time behavior of
Feynman-Kac semigroups, leading to a precise description of the asymptotic
behavior of the variance term $v_{n}$ appearing in Theorem \ref{lognormal} in
two special cases : time-homogeneous models, and models in a stationary
ergodic random environment.

The key result, Theorem (\ref{lognormal}), is established in Section
\ref{sec:fluctuations}. The key idea is to expand $\log{\overline{\gamma}%
_{n}^{N}(1)}$ in terms of local fluctuation terms of the form $V_{k}^{N}$.
Broadly speaking, the contribution of quadratic terms in the expansion amounts
to an asymptotically deterministic bias term whose fluctuations are controlled
with variance bounds, while the contribution of linear terms is treated by
invoking the martingale central limit theorem.

\section{Regularity of the covariance function}

\label{cov-regularity}

We first note that, in the special case where $K_{n,\eta}%
(x,\mbox{\LARGE .})=\Phi_{n}(\eta)$ for all $x$, Property \eqref{H2} is in
fact a consequence of \eqref{H0} and \eqref{H}. Indeed, we can then write
\[
\lbrack\Phi_{n}(\mu_{1})-\Phi_{n}(\mu_{2})](f)=\frac{1}{\mu_{1}(G_{n-1})}%
~[\mu_{1} -\mu_{2}]\left(  G_{n-1}M_{n}(f-\Phi_{n}(\mu_{2})(f))\right)
\]
and check that, for all $f \in\mathrm{Osc}(E)$, one has
\[
\Vert\left[  K_{n,\mu_{1}}-K_{n,\mu_{2}}\right]  (f)\Vert\leq2g~\left\vert
\left[  \mu_{1}-\mu_{2}\right]  (h_{n,\mu_{2}})\right\vert ,
\]
where $g$ is defined in (\ref{H}) and
\[
h_{n,\mu}=\frac{1}{2\Vert G_{n-1}\Vert}~G_{n-1}M_{n}(f-\Phi_{n}(\mu
)(f))\in\mbox{\rm
Osc}(E).
\]

%Note that we have used the following general property: for any two probability measures $\mu_1, \mu_2$ on $E$, one has that
%\begin{equation}\label{tv-osc}  \Vert \mu_1 - \mu_2  \Vert_{{\tiny \mathrm{tv}}}  = \sup_{f \in {\rm Osc}(E)}  \vert   \mu_1(f) - \mu_2(f)  \vert.\end{equation}
In the alternative case (\ref{2nd-choice}), we have
\[
\left[  K_{n,\mu_{1}}-K_{n,\mu_{2}}\right]  (f)=\left(  1-G_{n-1}\right)
~\left[  \Phi_{n}(\mu_{1})-\Phi_{n}(\mu_{2})\right]  (f)
\]
so that \eqref{H2} is also satisfied.

Observe that \eqref{H2} immediately implies the following Lipschitz-type
property:
\begin{equation}
\sup_{x \in E} \Vert K_{n,\mu_{1}}(x, \cdot) -K_{n,\mu_{2}}(x,\cdot)
\Vert_{{\tiny \mathrm{tv}}} \leq\kappa\Vert\mu_{1} - \mu_{2} \Vert
_{{\tiny \mathrm{tv}}} . \label{LipK}%
\end{equation}

\begin{prop}
\label{Lip-Cov} One has the following bound, valid for any two probability
measures $\mu_{1}, \mu_{2}$ on $E$, and functions $f_{1}, f_{2} \in
\mathrm{Osc}(E)$:
\begin{equation}
\label{formule-cov}\vert\mbox{\rm Cov}_{n,\mu_{1}}(f_{1},f_{2}%
)-\mbox{\rm Cov}_{n,\mu_{2}}(f_{1},f_{2}) \vert\leq c \Vert\mu_{1} - \mu_{2}
\Vert_{{\tiny \mathrm{tv}}}.
\end{equation}

\end{prop}

\noindent\mbox{\bf Proof:}\newline

We have
\[%
\begin{array}
[c]{l}%
\mbox{\rm Cov}_{n,\mu_{1}}(f_{1},f_{2})-\mbox{\rm Cov}_{n,\mu_{2}}(f_{1}%
,f_{2})\\
\\
=\left[  \Phi_{n}(\mu_{1})-\Phi_{n}(\mu_{2})\right]  (f_{1}f_{2})+[\mu_{2}%
-\mu_{1}]\left(  K_{n,\mu_{2}}(f_{1})K_{n,\mu_{2}}(f_{2})\right) \\
\\
\hskip5cm+\mu_{1}\left(  K_{n,\mu_{2}}(f_{1})K_{n,\mu_{2}}(f_{2})-K_{n,\mu
_{1}}(f_{1})K_{n,\mu_{1}}(f_{2})\right)
\end{array}
\]
and
\[
\lbrack\Phi_{n}(\mu_{1})-\Phi_{n}(\mu_{2})]=\mu_{1}\lbrack K_{n,\mu_{1}%
}-K_{n,\mu_{2}}]+[\mu_{1}-\mu_{2}]K_{n,\mu_{2}}.
\]
Note that there is no loss of generality in assuming that $\mu_{2}(f_{1}%
)=\mu_{2}(f_{2})=0$, so that $\Vert f_{i}\Vert\leq\mbox{\rm osc}(f_{i})\leq1$.
Thus, using
\[%
\begin{array}
[c]{l}%
\left\Vert K_{n,\mu_{2}}(f_{1})K_{n,\mu_{2}}(f_{2})-K_{n,\mu_{1}}%
(f_{1})K_{n,\mu_{1}}(f_{2})\right\Vert \\
\\
\leq\left\Vert K_{n,\mu_{1}}(f_{1})-K_{n,\mu_{2}}(f_{1})\right\Vert
+~\left\Vert K_{n,\mu_{1}}(f_{2})-K_{n,\mu_{2}}(f_{2})\right\Vert ,
\end{array}
\]
the desired conclusion follows from \eqref{LipK}.

\hfill\hbox{\vrule
height 5pt width 5pt depth 0pt}\medskip\newline

We also state the easily checked Lipschitz type bound, valid for all
$f_{1},f_{2}, \phi_{1}, \phi_{2} \in\mathcal{B}_{b}(E)$
\begin{equation}%
\begin{array}
[c]{l}%
\left\vert \mbox{\rm Cov}_{n,\mu}(f_{1},f_{2})-\mbox{\rm Cov}_{n,\mu}(\phi
_{1},\phi_{2})\right\vert \\
\\
\leq c \left(  \Vert f_{1}\Vert~\Vert f_{2}-\phi_{2}\Vert+\Vert\phi_{2}
\Vert~\Vert f_{1}-\phi_{1}\Vert\right)  .
\end{array}
\label{formule-cov-f}%
\end{equation}

\section{Feynman-Kac semigroups}

\label{sec-quantitative}

\subsection{Contraction estimates}

\label{sec-quantitative-1} We denote by $(\Phi_{p,n})_{0 \leq p \leq n}$ the
semigroup of nonlinear operators acting on probability measures defined by
\[
\Phi_{p,n} := \Phi_{n} \circ\cdots\circ\Phi_{p+1},
\]
so that
\begin{equation}
{\eta_{n}}\left(  f\right)  =\Phi_{p,n}(\eta_{p})(f)={\eta_{p}Q_{p,n}%
(f)}/{\eta_{p}Q_{p,n}(1)}=\Psi_{Q_{p,n}(1)}(\eta_{p})P_{p,n}(f). \label{Ppnb}%
\end{equation}
One has that
\begin{equation}
\sup_{\mu,\nu}{\Vert\Phi_{p,n}(\mu)-\Phi_{p,n}(\nu)\Vert_{{\tiny \mathrm{tv}}%
}}=\beta(P_{p,n}), \label{contract-Phi}%
\end{equation}
see for example \cite[chapter 4]{DMo04}. We also set
\[
g_{p,n}:=\sup_{x,y\in E}{\left[  Q_{p,n}(1)(x)/Q_{p,n}(1)(y)\right]  }%
\quad\mbox{\rm and}\quad d_{p,n}(f)=\overline{Q}_{p,n}(f-\eta_{n}(f)).
\]
Note that $\overline{Q}_{n,n+1}(1)=G_{n}/\eta_{n}(G_{n})=\overline{G}_{n}$,
and that
\begin{equation}
d_{p,n}(\overline{G}_{n})=\overline{Q}_{p,n}(\overline{Q}_{n,n+1}%
(1)-1)=\overline{Q}_{p,n+1}(1)-\overline{Q}_{p,n}(1). \label{Qpn-diff}%
\end{equation}
We will use the fact that the semigroup $Q_{p,n}$ satisfies a decomposition
similar to (\ref{multiplicative}): for any probability measure $\mu$ on $E$,
one has that
\begin{equation}
\mu Q_{p,n}(1)=\prod_{p\leq q<n}{\Phi_{p,q}(\mu)(G_{q}).} \label{eq:repQpn}%
\end{equation}

\label{sec-quantitative-2} Also, combining (\ref{eq:normalizedQpn}) and
(\ref{eq:repQpn}), we can write
\begin{equation}
\label{decomp-log}\log{\overline{Q}_{p,n}(1)(x)}=\sum_{p\leq q<n}\left[
\log{\Phi_{p,q}(\delta_{x})(G_{q})}-\log{\Phi_{p,q}(\eta_{p})(G_{q})}\right]
.
\end{equation}

\begin{lem}
\label{lem-cv-exp} For any $0\leq p\leq n$ and any $f\in\mbox{\rm Osc}(E)$, we
have
\begin{equation}
g_{p,n}\leq b:=\exp{\left(  a(g-1)/(1-e^{-\lambda})\right)  }\quad
\mbox{and}\quad\left\Vert d_{p,n}(f)\right\Vert \leq ab~e^{-\lambda(n-p)}.
\label{control-gpn}%
\end{equation}
In addition, for any $\mu,\nu\in\mathcal{P}(E)$ we have
\begin{equation}
\Vert\Phi_{p,n}(\mu)-\Phi_{p,n}(\nu)\Vert_{{\tiny \mathrm{tv}}}\leq
ab~e^{-\lambda(n-p)}~\Vert\mu-\nu\Vert_{{\tiny \mathrm{tv}}}.
\label{contract-Lipschitz-Phi}%
\end{equation}

\end{lem}

\noindent\mbox{\bf Proof:}\newline Using the decomposition (\ref{eq:repQpn}),
we have
\begin{equation}
\frac{Q_{p,n}(1)(x)}{Q_{p,n}(1)(y)}=\frac{{\delta_{x}}Q_{p,n}(1)}{{\delta_{y}%
}Q_{p,n}(1)}=\exp{\left\{  \sum_{p\leq q<n}\left(  \log{\Phi_{p,q}(\delta
_{x})(G_{q})}-\log{\Phi_{p,q}(\delta_{y})(G_{q})}\right)  \right\}  .}
\label{Qpnratio}%
\end{equation}
From the identity $\log{u}-\log{v}=\int_{0}^{1}\frac{(u-v)}{u+t(v-u)}~dt$,
valid for any $u,v>0$, we deduce the inequality
\[%
\begin{array}
[c]{l}%
\displaystyle\frac{Q_{p,n}(1)(x)}{Q_{p,n}(1)(y)}\leq\displaystyle\exp{\left\{
\sum_{p\leq q<n}~\displaystyle\widetilde{g}_{q}\times\left|  \Phi_{p,q}%
(\delta_{x})(\widetilde{G}_{q})-\Phi_{p,q}(\delta_{y})(\widetilde{G}%
_{q})\right|  \right\}  },
\end{array}
\]
with $\widetilde{G}_{q}:=G_{q}/\mbox{\rm osc}(G_{q})$ (and the convention that
$\widetilde{G}_{q}:=1$ if $G_{q}$ is constant), and $\widetilde{g}%
_{q}:={\mbox{\rm osc}(G_{q})}/{\inf G_{q}}\leq g_{q}-1$.

Using \eqref{H} and \eqref{contract-Phi}, we deduce that
\[
g_{p,n}\leq\exp{\left\{  a(g-1)\sum_{p\leq q<n}~e^{-\lambda(q-p)}\right\}
}\leq b.
\]
This ends the proof of the l.h.s. of (\ref{control-gpn}). The proof of the
r.h.s. of (\ref{control-gpn}) comes from the following expression for
$d_{p,n}(f)$:
\[
d_{p,n}(f)=\overline{Q}_{p,n}(1)\times P_{p,n}\left[  f-\Psi_{Q_{p,n}(1)}%
(\eta_{p})P_{p,n}(f)\right]
\]
which implies, using the fact that $\Vert\overline{Q}_{p,n}(1) \Vert\leq
g_{p,n}$, that
\begin{equation}
\left\Vert d_{p,n}(f)\right\Vert \leq g_{p,n}~\beta(P_{p,n}%
)~\mbox{\rm osc}(f)\leq ab~e^{-\lambda(n-p)}~\mbox{\rm osc}(f).
\label{dpn-estimate}%
\end{equation}
From \cite[Section 4.3]{DMo04}, see also Proposition 3.1 in~\cite{dhw-2012},
we have
\[
\Vert\Phi_{p,n}(\mu)-\Phi_{p,n}(\nu)\Vert_{{\tiny \mathrm{tv}}}\leq
g_{p,n}~\beta(P_{p,n})~\Vert\mu-\nu\Vert_{{\tiny \mathrm{tv}}}.
\]
Using (\ref{control-gpn}), we conclude that
\[
\Vert\Phi_{p,n}(\mu)-\Phi_{p,n}(\nu)\Vert_{{\tiny \mathrm{tv}}}\leq
ab~e^{-\lambda(n-p)}~\Vert\mu-\nu\Vert_{{\tiny \mathrm{tv}}}.
\]
This ends the proof of the lemma. \hfill\hbox{\vrule height 5pt width 5pt
depth 0pt}\medskip\newline

\subsection{Limiting semigroup}

We now state a general theorem on the convergence of $\overline{Q}_{p,n}(1)$
when $n \to+\infty$.

\begin{theo}
\label{theo-Qpn} The following bound holds for all $0\leq p\leq n$:
\begin{equation}
\left\Vert \overline{Q}_{p,n}(1)-\overline{Q}_{p,\infty}(1)\right\Vert \leq
c~e^{-\lambda(n-p)}, \label{ref-limQpn}%
\end{equation}
where the limiting function $\overline{Q}_{p,\infty}(1)$ is defined through
the following series:
\begin{equation}
\log{\overline{Q}_{p,\infty}(1)(x)}:=\sum_{q\geq p}\left[  \log{\Phi
_{p,q}(\delta_{x})(G_{q})}-\log{\Phi_{p,q}(\eta_{p})(G_{q})}\right]  .
\label{def-Qinfty}%
\end{equation}

\end{theo}

\noindent\textbf{Proof of Theorem~\ref{theo-Qpn}:}\newline We first check that
the function $\overline{Q}_{p,\infty}(1)$ is well defined, using the fact
that, as in the proof of Lemma \ref{lem-cv-exp},
\[
\left|  \log{\Phi_{p,q}(\delta_{x})(G_{q})}-\log{\Phi_{p,q}(\eta_{p})(G_{q})}
\right|  \leq a(g-1) e^{-\lambda(q-p)}.
\]
One then has that
\[
\left|  \log{\ \overline{Q}_{p,n}(1)}(x)-\log{\ \overline{Q}_{p,\infty}
(1)(x)}\right|  \leq\sum_{q\geq n} \left|  \log{\Phi_{p,q}(\delta_{x})(G_{q}%
)}-\log{\Phi_{p,q}(\eta_{p})(G_{q})} \right|  ,
\]
whence
\[
\left|  \log{\ \overline{Q}_{p,n}(1)}(x)-\log{\ \overline{Q}_{p,\infty}
(1)(x)}\right|  \leq\sum_{q\geq n} a (g-1)~e^{-\lambda(q-p)} \leq
c~e^{-\lambda(n-p)}.
\]

Using the identity $e^{u}-e^{v}=(x-y)~\int_{0}^{1}e^{tu+(1-t)v}dt$, we finally
check that
\[
\left\Vert \overline{Q}_{p,n}(1)-\overline{Q}_{p,\infty}(1)\right\Vert \leq
c~\left\Vert \log{\ \overline{Q}_{p,n}(1)}-\log{\ \overline{Q}_{p,\infty}
(1)}\right\Vert ,
\]
thanks to the fact that $\Vert\overline{Q}_{p,n}(1) \Vert\leq g_{p,n} \leq g$.
This ends the proof of (\ref{ref-limQpn}). \hfill\hbox{\vrule
height 5pt width 5pt depth 0pt}\medskip\newline

\subsection{The time-homogeneous case}

\label{time-homog}

Here we consider the special case of time-homogeneous models, where there
exist $G, M, K$ such that $G_{n}=G$ for all $n \geq0$, and $M_{n}=M$ and
$K_{n}=K$ for all $n \geq1$.

Our assumptions imply the existence of a unique fixed point $\eta_{\infty
}=\Phi(\eta_{\infty})$ towards which $\eta_{n}$ converges exponentially fast:
for all $n \geq0$,
\begin{equation}
\Vert\Phi^{n}(\eta_{0})-\eta_{\infty}\Vert_{{\tiny \mathrm{tv}}}\leq
ab~e^{-\lambda n}. \label{fixed-point}%
\end{equation}

In this situation, Theorem \ref{theo-Qpn} leads to a precise description of
the asymptotic behavior of the variance term $v_{n}$ appearing in Theorem
\ref{lognormal}. To state it, consider the fixed point measure $\eta_{\infty}$
introduced in (\ref{fixed-point}), and define the function $h$ by
\[
\log{h(x)}:=\sum_{n \geq0}\left[  \log{\Phi^{n}(\delta_{x})(G_{n})}-\log
{\Phi^{n}(\eta_{\infty})(G_{q})}\right]  .
\]
In the stationary version of the model where $\eta_{0}:=\eta_{\infty}$, $h$
corresponds to the limiting function $\overline{Q}_{0,\infty}(1)$ whose
existence is asserted by Theorem \ref{theo-Qpn}. In this situation, it turns
out that, by stationarity, $\overline{Q}_{n,\infty}(1) = h$ for all $n \geq1$.

\begin{prop}
\label{p:cv-Qinfini-h} One has the following bound for all $p \geq0$:
\begin{equation}
\label{remplacement-Qh}|| \overline{Q}_{p,\infty}(1) - h || \leq c e^{-\lambda
p}.
\end{equation}

\end{prop}

\begin{cor}
\label{corollary:Qlimit} One has that
\begin{equation}
\frac{1}{n}\sum_{0\leq p<n}\mbox{\rm Cov}_{p,\eta_{p-1}}(\overline{Q}
_{p,n}(1),\overline{Q}_{p,n}(1))=\mbox{\rm Cov}_{\eta_{\infty}}
(h,h)+\mbox{\rm O}(1/n) , \label{ref-limi-Cov}%
\end{equation}
where we use the notation $\mbox{\rm Cov}_{\eta}$ to denote the common value
of $\mbox{\rm Cov}_{p,\eta}$ for $p \geq1$.
\end{cor}

An alternative spectral characterization of the map $h$ is given in the
following corollary. In the homogeneous case, $Q_{p,p+1}$ does not depend on
$p$, so we use the simpler notation $Q$.

\begin{cor}
\label{corollary:spectral} In the homogeneous case, the $(\eta_{\infty} Q(1),
h)$ is characterized as the unique pair $(\zeta, f)$ such that $Q(f) = \zeta
f$ and $\eta_{\infty}(f)=1$.
\end{cor}

\noindent\textbf{Proof of Proposition~\ref{p:cv-Qinfini-h}:}\newline Using the
exponential convergence to $\eta_{\infty}$ stated in \eqref{fixed-point}, and
the Lipschitz property \eqref{contract-Lipschitz-Phi}, we have that
\[
\sum_{q\geq p}\left[  \log{\Phi_{p,q}(\eta_{p})(G_{q})}-\log{\Phi_{p,q}%
(\eta_{\infty})(G_{q})}\right]  \leq c~\sum_{q\geq p}~e^{-\lambda
((q-p)+p)}\leq c^{\prime}~e^{-\lambda p}.
\]
We conclude as in the proof of Theorem \ref{theo-Qpn}. \hfill\hbox{\vrule
height 5pt width 5pt depth 0pt}\medskip\newline

\noindent\textbf{Proof of Corollary~\ref{corollary:Qlimit}:}\newline

Using the Lipschitz property \eqref{formule-cov}, and the fact that, for all
$p,n$, $\Vert\overline{Q}_{p,n}(1) \Vert\leq g$, we see that replacing each
$\eta_{p-1}$ in the l.h.s. of \eqref{ref-limi-Cov} by $\eta_{\infty}$ leads to
a $O(1/n)$ error term. Then, using Theorem \ref{theo-Qpn} and
\eqref{formule-cov-f}, we see that we can replace each $\overline{Q}_{p,n}(1)$
term by $\overline{Q}_{p,\infty}(1)$ in the l.h.s. of \eqref{ref-limi-Cov},
and commit no more than a $O(1/n)$ overall error. Finally,
\eqref{remplacement-Qh}, allows us to replace each $\overline{Q}_{p,\infty
}(1)$ by $h$, again with an overall $O(1/n)$ error term. \hfill\hbox{\vrule
height 5pt width 5pt depth 0pt}\medskip\newline

\noindent\textbf{Proof of Corollary~\ref{corollary:spectral}:}\newline

We consider the stationary version of the model where we start with $\eta
_{0}:=\eta_{\infty}$.

Let us first check that one indeed has $\eta_{\infty}(h)=1$ and $Q(h) =
\eta_{\infty}(Q(1)) h$. By Theorem \ref{theo-Qpn}, we have that
\begin{equation}
\lim_{n \to+\infty} \left\Vert \overline{Q}_{0,n}(1) - h \right\Vert = 0.
\label{cv-h}%
\end{equation}
Since by construction, $\eta_{\infty} \overline{Q}_{0,n}(1) =1$, \eqref{cv-h}
yields that $\eta_{\infty}(h)=1$. Then, due to stationarity, one has
$\overline{Q}_{p,n} = \overline{Q}^{n-p}$, with $\overline{Q}(f) := Q(f) /
\eta_{\infty} Q(1)$, so that one can also deduce from \eqref{cv-h} that
$\overline{Q}(h) = h$, which yields that $Q(h) = \eta_{\infty}(Q(1)) h$.

Now consider a pair $(\zeta, f)$ such that $Q(f) = \zeta f$ and $\eta_{\infty
}(f)=1$, and let us show that $\zeta= \eta_{\infty} Q(1)$ and $f=h$.

By stationarity, one has that
\[
\overline{Q}_{0,n}(f) = Q^{n}(f) / \eta_{\infty} Q^{n}(1),
\]
and we deduce from \eqref{eq:repQpn} and the stationarity of $\eta_{\infty}$
that
\[
\eta_{\infty} Q^{n}(1) = \left(  \eta_{\infty} Q(1)\right)  ^{n}.
\]

Using the fact that $\Phi(\eta_{\infty}) = \eta_{\infty}$, we have the
identity
\[
\eta_{\infty} Q (f) / \eta_{\infty}(Q(1)) = \eta_{\infty}(f).
\]
Since $Q(f) = \zeta f$ and $\eta_{\infty}(f)=1$, we immediately deduce that
$\zeta= \eta_{\infty}(Q(1))$.

As a consequence, the fact that $Q(f) = \eta_{\infty}(Q(1)) f$ implies that,
for all $n \geq1$, one has
\[
\overline{Q}_{0,n}(f) = f.
\]

On the other hand, given two bounded functions $f_{1},f_{2}$, we have that
\[
\overline{Q}_{0,n}(f_{1}-f_{2})(x)=\Phi_{0,n}(\delta_{x})(f_{1}-f_{2}%
)\times\overline{Q}_{0,n}(1)(x)
\]
Letting $n\rightarrow\infty$, \eqref{fixed-point} and Theorem \ref{theo-Qpn}
yield
\[
\lim_{n\rightarrow\infty}\overline{Q}_{0,n}(f_{1}-f_{2})=\eta_{\infty}%
(f_{1}-f_{2})\times h.
\]
Using $f_{1}:=f$ and $f_{2}:=h$, we deduce that $f=h$. \hfill\hbox{\vrule
height 5pt width 5pt depth 0pt}\medskip\newline

\subsection{The random environment case}

\label{FK-random-environment}

\subsubsection{Description of the model}

We consider a stationary and ergodic process $Y=(Y_{n})_{n\in\mathbb{Z}}$
taking values in a measurable state space $(S, \mathcal{S})$. The process $Y$
provides a random environment governing the successive transitions between
step $n-1$ and step $n$ in our model. In the sequel, we define and study the
model for a given realization $y \in S^{\mathbb{Z}}$ of the environment. It is
only in Corollary \ref{limi-Cov-re} that we exploit the ergodicity of $Y$ to
establish the almost sure limiting behavior of the variance $v_{n}$.

Specifically, we consider a family $(M_{s})_{s\in S}$ of Markov kernels on
$E$, a family $(G_{s})_{s \in S}$ of positive bounded functions on $E$.

For $n \in\mathbb{Z}$ and $y \in S^{\mathbb{Z}}$, we set $M_{n}^{y} :=
M_{y_{n}}$ and $G_{n}^{y} := G_{y_{n}}$. We then denote with a $y$ superscript
all the objects associated with the Feynman-Kac model using the sequence of
kernels $(M_{n}^{y})_{n \geq1}$ and functions $(G_{n}^{y})_{n \geq0}$, i.e.
the measures $\gamma_{n}^{y}$ and $\eta_{n}^{y}$, the operators $\Phi
^{y}_{p,n}$, $G^{y}_{p,n}$, $\mathrm{Cov}^{y}_{p,\eta}$, etc. To define the
particle approximation scheme, we also consider a family of Markov kernels
$(K_{(s,s^{\prime}), \mu})_{s,s^{\prime}\in S, \ \mu\in\mathcal{P}(E)}$ such
that, for all $s,s^{\prime},\mu$, one has
\[
\Psi_{G_{s}}(\mu) M_{s^{\prime}} = \mu K_{(s,s^{\prime}), \mu}.
\]
We then use $K^{y}_{n, \mu} := K_{(y_{n-1},y_{n}), \mu}$ for all $n \geq1$.

We then define the shift operator on $S^{\mathbb{Z}}$ by setting, for every $y
= (y_{n})_{n \in\mathbb{Z}} \in S^{\mathbb{Z}}$, $\theta(y) := (y_{n+1})_{n
\in\mathbb{Z}}$. With our definitions, one has that, for all $0 \leq p \leq
n$,
\[
Q^{y}_{p,n} = Q_{0,n-p}^{\theta^{p}(y)}, \ \Phi^{y}_{p,n} = \Phi
_{0,n-p}^{\theta^{p}(y)},
\]
and in particular
\begin{equation}
\label{composhift}\Phi^{y}_{0,n} = \Phi^{\theta^{p}(y)}_{0,n-p} \circ\Phi
^{y}_{0,p}.
\end{equation}

Our assumptions on the model are that $E$ has a Polish space structure, and
that the bounds listed in \eqref{H0}, \eqref{H} and \eqref{H2} hold for
$M_{n}^{y}$, $G_{n}^{y}$ and $K_{n,\mu}^{y}$ uniformly over $y \in
S^{\mathbb{Z}}$.

\subsubsection{Contraction properties}

Rewriting (\ref{contract-Phi}) and (\ref{contract-Lipschitz-Phi}) in the
present context, we have that, for all $y$,
\begin{equation}
\beta\left(  P_{0,n}^{y}\right)  =\sup_{\mu,\nu}{\Vert\Phi_{0,n}^{y}(\mu
)-\Phi_{0,n}^{y}(\nu)\Vert_{{\tiny \mathrm{tv}}}}\leq a~e^{-\lambda n}
\label{r0}%
\end{equation}
and
\begin{equation}
\label{contract-Lipschitz-Phi-y}{\Vert\Phi_{0,n}^{y}(\mu)-\Phi_{0,n}^{y}%
(\nu)\Vert_{{\tiny \mathrm{tv}}}}\leq ab~e^{-\lambda n}~\Vert\mu-\nu
\Vert_{{\tiny \mathrm{tv}}},
\end{equation}
with the constant $b$ defined in (\ref{control-gpn}). Using (\ref{composhift}%
), we have
\[
\Phi_{0,n+m}^{\theta^{-(n+m)}(y)}=\Phi_{0,n}^{\theta^{-n}(y)}\circ\Phi
_{0,m}^{\theta^{-(n+m)}(y)},
\]
so that, using \eqref{r0}, one has that
\[
\sup_{\mu,\nu}{\Vert\Phi_{0,n}^{\theta^{-n}(y)}(\mu)-\Phi_{0,n+m}%
^{\theta^{-(n+m)}(y)}(\nu)\Vert_{{\tiny \mathrm{tv}}}}\leq a~e^{-\lambda n}.
\]
Arguing as in~\cite{grey,nick}, we conclude that for any $f\in\mathcal{B}%
_{b}(E)$, and any $\mu\in\mathcal{P}(E)$, $\Phi_{0,n}^{\theta^{-n}(y)}%
(\mu)(f)$ is a Cauchy sequence, so that $\boldsymbol{\Phi}_{0,n}^{\theta
^{-n}(y)}(\mu)$ weakly converges to a measure $\eta_{\infty}^{y}$, as $n
\to\infty$. In addition, for any $n\geq0$, we have
\begin{equation}
\Phi_{0,n}^{y}(\eta_{\infty}^{y})=\eta_{\infty}^{\theta^{n}(y)}
\label{fixed-point-Y}%
\end{equation}
and exponential convergence to equilibrium
\begin{equation}
\label{exp-cv-eq-y}\sup_{\mu}{\Vert\Phi_{0,n}^{\theta^{-n}(y)}(\mu
)-\eta_{\infty}^{y}\Vert_{{\tiny \mathrm{tv}}}}\leq a~e^{-\lambda n}.
\end{equation}

We now restate the conclusion of Theorem~\ref{theo-Qpn} in the present context
: for all $0\leq p\leq n$, one has that
\begin{equation}
\left\Vert \overline{Q}^{y}_{p,n}(1)-\overline{Q}^{y}_{p,\infty}(1)\right\Vert
\leq c~e^{-\lambda(n-p)}, \label{ref-limQpn-y}%
\end{equation}
where the limiting function $\overline{Q}^{y}_{p,\infty}(1)$ is defined
through the series:
\begin{equation}
\log{\overline{Q}_{p,\infty}^{y}(1)(x)}:=\sum_{q\geq p}\left[  \log{\Phi^{y}
_{p,q}(\delta_{x})(G^{y}_{q})}-\log{\Phi^{y}_{p,q}(\eta_{p}^{y})(G_{q}^{y}%
)}\right]  . \label{def-Qinfty-y}%
\end{equation}

We now define the map $h^{y}$ by
\[
h^{y}(x):= \sum_{q \geq0}\left[  \log{\Phi^{y} _{0,q}(\delta_{x})(G^{y}_{q}%
)}-\log{\Phi^{y}_{0,q}(\eta_{\infty}^{y})(G_{q}^{y})}\right]  .
\]

\begin{prop}
\label{p:cv-Qinfini-hy} One has the following bound, valid for all $y \in
S^{\mathbb{Z}}$ and $p \geq0$:
\begin{equation}
\label{remplacement-Qh-re}|| \overline{Q}_{p,\infty}^{y}(1) - h^{\theta
^{p}(y)} || \leq c e^{-\lambda p}.
\end{equation}

\end{prop}

\textbf{Proof of Proposition~\ref{p:cv-Qinfini-hy}:}\newline

Setting $q:=q-p$ in the definition, we rewrite
\[
\log{\overline{Q}_{p,\infty}^{y}(1)(x)}=\sum_{q \geq0}\left[  \log
{\Phi^{\theta^{p}(y)} _{0,q}(\delta_{x})(G^{\theta^{p}(y)}_{q})}-\log
{\Phi^{\theta^{p}(y)}_{0,q}(\eta_{p}^{y})(G_{q}^{\theta^{p}(y)})}\right]  .
\]
On the other hand,
\[
h^{\theta^{p}(y)}(x)= \sum_{q \geq0}\left[  \log{\Phi^{\theta^{p}(y)}
_{0,q}(\delta_{x})(G^{\theta^{p}(y)}_{q})}-\log{\Phi^{\theta^{p}(y)}%
_{0,q}(\eta_{\infty}^{\theta^{p}(y)})(G_{q}^{\theta^{p}(y)})}\right]  .
\]

Using \eqref{exp-cv-eq-y}, we obtain that
\[
\Vert\eta_{p}^{y} -\eta_{\infty}^{\theta^{p}(y)}\Vert_{{\tiny \mathrm{tv}}%
}\leq a~e^{-\lambda p}.
\]
Combining this bound with \eqref{contract-Lipschitz-Phi-y}, we deduce that
\[
\left|  \Phi^{\theta^{p}(y)}_{0,q}(\eta_{p}^{y})(G_{q}^{\theta^{p}(y)}) -
\Phi^{\theta^{p}(y)}_{0,q}(\eta_{\infty}^{\theta^{p}(y)})(G_{q}^{\theta
^{p}(y)}) \right|  \leq c e^{-\lambda(p+q)}.
\]
We conclude as in the proof of Theorem \ref{theo-Qpn}. \hfill
\hbox{\vrule height 5pt width 5pt depth 0pt}\medskip\newline

Introduce the map $\mathcal{C}$ defined on $S^{\mathbb{Z}}$ by
\[
\mathcal{C}(y) := \mbox{\rm Cov}^{\theta^{-1}(y)}_{1,\eta_{\infty}%
^{\theta^{-1}(y)} }(h^{y},h^{y}).
\]
We add to \eqref{H0}-\eqref{H}-\eqref{H2} the assumption that $\mathcal{C}$ is
measurable with respect to the product $\sigma-$algebra on $S^{\mathbb{Z}}$.

Arguing as in the proof of Corollary \ref{ref-limi-Cov}, then applying the
ergodic theorem, we deduce the following asymptotic behavior for the variance
$v_{n}$.

\begin{cor}
\label{limi-Cov-re} One has the following bound:
\[%
\begin{array}
[c]{l}%
\frac{1}{n}\sum_{0\leq p<n}\mbox{\rm Cov}^{y}_{p,\eta_{p-1}^{y}}(\overline
{Q}_{p,n}^{y}(1),\overline{Q}_{p,n}^{y}(1))\\
\\
=\frac{1}{n}\sum_{1\leq p<n} \mathcal{C}(\theta^{p}(y)) + O(1/n).
\end{array}
\]

In addition, we have
\[%
\begin{array}
[c]{l}%
\displaystyle\lim_{n\rightarrow\infty}\frac{1}{n}\sum_{0\leq p<n}\mbox{\rm
Cov}_{p,\eta_{p-1}^{Y}}(\overline{Q}_{p,n}^{Y}(1),\overline{Q}_{p,n}%
^{Y}(1)))=\mathbb{E}\left(  \mbox{\rm Cov}^{\theta^{-1}(Y)}_{1,\eta_{\infty
}^{\theta^{-1}(Y)}}(h^{Y},h^{Y})\right)  a.s.
\end{array}
\]

\end{cor}

\section{Fluctuation analysis\label{sec:fluctuations}}

\subsection{Moment bounds}

In addition to the local error fields $V^{N}_{n}$ defined in
\eqref{VnNlocalsamplingerror}, we consider the global error fields $W^{N}_{n}$
defined by
\begin{equation}
W^{N}_{n}=\sqrt{N}\text{ }\left(  \eta_{n}^{N}-\eta_{n}\right)  \text{
}\Leftrightarrow\text{ }\eta_{n}^{N}=\eta_{n}+\frac{1}{\sqrt{N}}~W^{N}_{n}.
\label{Wnerror}%
\end{equation}

We now quote key moment estimates on $V^{N}_{n}$ and $W^{N}_{n}$, see
\cite[chapter 4]{DMo04} or \cite[chapter 9]{dp-2013}. Under our assumptions,
one has that, for all $n \geq0$, $N \geq1$, all $f\in\mbox{\rm Osc}(E)$ and $m
\geq1$,
\begin{equation}
\left\Vert V_{n}^{N}(f)\right\Vert _{m}\leq c(m), \label{Lm-estimate}%
\end{equation}
and
\begin{equation}
\label{W-estimate}\left\Vert W_{n}^{N}(f)\right\Vert _{m}\leq c(m).
\end{equation}

\subsection{Expansion of the particle estimate of log-normalizing
constants\label{sec:expansionfirst}}

Starting from the product-form expression (\ref{relativeestimate}), we apply a
second-order expansion for the logarithm of each factor. Using
\eqref{W-estimate}, we have that, for all $n \geq0$ and $N \geq1$,
\begin{align}
\log\overline{\gamma}_{n}^{N}(1)  &  =\frac{1}{\sqrt{N}}\sum_{0\leq p<n}%
W_{p}^{N}(\overline{G}_{p})-\frac{1}{2N}\sum_{0\leq p<n}\left(  W_{p}%
^{N}(\overline{G}_{p})\right)  ^{2}\nonumber\\
&  +\frac{1}{\sqrt{N}}\left(  \frac{n}{N}\right)  C(n,N), \label{develop-log}%
\end{align}
where, for all $m \geq1$, the remainder term satisfies the moment
$||C(n,N)||_{m}\leq c(m)$.

\subsection{Second order perturbation formulae\label{decomp-de-W}}

We derive an expansion of $W_{n}^{N}(f)$ in terms of local error terms
$V_{p}^{N}$ introduced in (\ref{VnNlocalsamplingerror}), up to an error term
of order $\frac{1}{N}$. The key result we prove is the following.

\begin{theo}
\label{theo-key} For all $n\geq0$, $N\geq1$ and any function $f\in
\mbox{\rm Osc}(E)$,
\begin{equation}%
\begin{array}
[c]{l}%
W_{n}^{N}(f)=W_{n}^{N}(f)+\displaystyle\frac{1}{N}~R_{n}^{N}(f)
\end{array}
\label{2dec}%
\end{equation}
where
\begin{align*}
W_{n}^{N}(f)  &  =\displaystyle\sum_{p=0}^{n}V_{p}^{N}\left[  d_{p,n}%
(f)\right] \\
&  -\displaystyle\frac{1}{\sqrt{N}}\sum_{0\leq p<n}~\left[  \sum_{q=0}%
^{p}V_{q}^{N}\left[  d_{q,p}(\overline{G}_{p})\right]  \right]  \left[
\sum_{q=0}^{p}V_{q}^{N}\left[  d_{q,n}(f)\right]  \right]
\end{align*}
and where the remainder measure $R_{n}^{N}$ is such that, for all $m\geq1$,
\[
||R_{n}^{N}(f)||_{m}\leq c(m).
\]

\end{theo}

To prove Theorem \ref{theo-key}, we start with the following exact
decomposition of $W_{n}^{N}(f)$ into a first term of order $1$ involving the
$V_{p}^{N}$ for $p=0,...,n$ plus a remainder term of order $1/\sqrt{N}$.

\begin{theo}
[{\cite[chapter 9]{dp-2013}}]\label{decomp-exacte} For all $n\geq0$, $N\geq1$
and any function $f\in\mbox{\rm Osc}(E)$, we have the decomposition
\begin{equation}%
\begin{array}
[c]{l}%
W_{n}^{N}(f)=\displaystyle\sum_{p=0}^{n}V_{p}^{N}\left[  d_{p,n}(f)\right]
+\displaystyle\frac{1}{\sqrt{N}}~S_{n}^{N}(f),
\end{array}
\label{heq-ok-form}%
\end{equation}
with the second order remainder
\[
S_{n}^{N}(f):=-\sum_{0\leq p<n}~\displaystyle\frac{1}{\eta_{p}^{N}%
(\overline{G}_{p})}~W_{p}^{N}(\overline{G}_{p})~W_{p}^{N}\left[
d_{p,n}(f)\right]  .
\]

\end{theo}

Note that, under our assumptions, the remainder term satisfies for all
$m\geq1$
\begin{equation}
||S_{n}^{N}(f)||_{m}\leq c(m). \label{premier-develop-W}%
\end{equation}
Decomposing $1/\eta_{p}^{N}(\overline{G}_{p})$ into a term of order $1$ plus a
term of order $1/\sqrt{N}$ as follows
\begin{equation}
\displaystyle\frac{1}{\eta_{p}^{N}(\overline{G}_{p})}=1-\displaystyle\frac
{1}{\eta_{p}^{N}(\overline{G}_{p})}\frac{1}{\sqrt{N}}~W_{p}^{N}(\overline
{G}_{p}), \label{dec-frac-etaQ}%
\end{equation}
we refine Theorem \ref{decomp-exacte} into the following decomposition, which
now has an error term of order $1/N$.

\begin{cor}
For all $n\geq0$, $N\geq1$ and any function $f\in\mbox{\rm Osc}(E)$, we have
the decomposition
\begin{equation}%
\begin{array}
[c]{l}%
W_{n}^{N}(f)\\
\\
=\displaystyle\sum_{p=0}^{n}V_{p}^{N}\left[  d_{p,n}(f)\right]
-\displaystyle\frac{1}{\sqrt{N}}\sum_{0\leq p<n}~W_{p}^{N}(\overline{G}%
_{p})~W_{p}^{N}\left[  d_{p,n}(f)\right]  +\frac{1}{N}\mathcal{R}_{n}^{N}(f)
\end{array}
\label{1dec}%
\end{equation}
where the remainder term is such that, for all $m\geq1$, $||\mathcal{R}%
_{n}^{N}(f)||_{m}\leq c(m)$.
\end{cor}

\noindent\mbox{\bf Proof:}\newline Using (\ref{dec-frac-etaQ}), we obtain
(\ref{1dec}) with the remainder term
\[
\mathcal{R}_{n}^{N}(f):=\sum_{0\leq p<n}~\displaystyle\frac{1}{\eta_{p}%
^{N}(\overline{G}_{p})}~W_{p}^{N}(\overline{G}_{p})^{2}~W_{p}^{N}\left[
d_{p,n}(f)\right]  .
\]
Note that, for any $m\geq1$, we have that
\[%
\begin{array}
[c]{l}%
\mathbb{E}\left(  \left\vert \mathcal{R}_{n}^{N}(f)\right\vert ^{m}\right)
^{\frac{1}{m}}\\
\\
\leq g~\sum_{0\leq p<n}\mathbb{E}\left(  \left\vert W_{p}^{N}(\overline{G}%
_{p})\right\vert ^{4m}\right)  ^{\frac{1}{2m}}\times\mathbb{E}\left(
\left\vert W_{p}^{N}\left[  d_{p,n}(f)\right]  \right\vert ^{2m}\right)
^{\frac{1}{2m}}.
\end{array}
\]
Combining (\ref{W-estimate}) and (\ref{control-gpn}), we find that
\[
\mathbb{E}\left(  \left\vert \mathcal{R}_{n}^{N}(f)\right\vert ^{m}\right)
^{\frac{1}{m}}\leq c\text{ }\sum_{0\leq p<n}e^{\lambda(n-p)}.
\]
This ends the proof of the corollary.\hfill\hbox{\vrule height 5pt width 5pt
depth 0pt}\medskip\newline

We are now ready to derive Theorem \ref{theo-key}, by replacing the $W^{N}%
_{p}$ terms appearing in the previous corollary by their expansions in terms
of the $V^{N}_{p}$ provided by Theorem \ref{decomp-exacte}. Here is the proof
of Theorem \ref{theo-key}.

\noindent\mbox{\bf Proof:}\newline Using (\ref{1dec}), we have
\[
W_{n}^{N}(f)=\mathcal{V}_{n}^{N}(f)+\displaystyle\frac{1}{\sqrt{N}%
}~\mathcal{W}_{n}^{N}(f)+\frac{1}{N}\mathcal{R}_{n}^{N}(f)
\]
with
\begin{align*}
\mathcal{V}_{n}^{N}(f):=  &  \displaystyle\sum_{p=0}^{n}V_{p}^{N}\left[
d_{p,n}(f)\right] \\
\mathcal{W}_{n}^{N}(f):=  &  -\sum_{0\leq p<n}~W_{p}^{N}(\overline{G}%
_{p})~W_{p}^{N}\left[  d_{p,n}(f)\right]
\end{align*}
This implies that
\[
\sum_{0\leq p<n}~W_{p}^{N}(\overline{G}_{p})~W_{p}^{N}\left[  d_{p,n}%
(f)\right]  =\mathcal{I}_{n}^{(0)}+\frac{1}{\sqrt{N}}~\mathcal{I}_{n}%
^{(1)}(f)+\frac{1}{N}~\mathcal{I}_{n}^{(2)}(f)+\frac{1}{N^{2}}~\mathcal{I}%
_{n}^{(3)}(f)
\]
with
\begin{align*}
\mathcal{I}_{n}^{(0)}(f)  &  =\sum_{0\leq p<n}\mathcal{V}_{p}^{N}(\overline
{G}_{p})~\mathcal{V}_{p}^{N}(d_{p,n}(f)),\\
\mathcal{I}_{n}^{(1)}(f)  &  =\sum_{0\leq p<n}~\left[  \mathcal{V}_{p}%
^{N}(\overline{G}_{p})~\mathcal{W}_{p}^{N}(d_{p,n}(f))+\mathcal{W}_{p}%
^{N}(\overline{G}_{p})~\mathcal{V}_{p}^{N}(d_{p,n}(f))\right]  ,\\
\mathcal{I}_{n}^{(2)}(f)  &  =\sum_{0\leq p<n}\left\{  \mathcal{R}_{p}%
^{N}(\overline{G}_{p})~\left[  \mathcal{V}_{p}^{N}(d_{p,n}%
(f))+\displaystyle\frac{1}{\sqrt{N}}~\mathcal{W}_{p}^{N}(d_{p,n}(f))\right]
\right. \\
&  \hskip3cm\left.  +\mathcal{R}_{p}^{N}(d_{p,n}(f))~\left[  \mathcal{V}%
_{p}^{N}(\overline{G}_{p})+\displaystyle\frac{1}{\sqrt{N}}~\mathcal{W}_{p}%
^{N}(\overline{G}_{p})\right]  \right\}  ,\\
\mathcal{I}_{n}^{(3)}(f)  &  =\sum_{0\leq p<n}\mathcal{R}_{p}^{N}(\overline
{G}_{p})~\mathcal{R}_{p}^{N}(d_{p,n}(f)).
\end{align*}
Arguing as in the previous proof, we see that $\sup_{1\leq i\leq3}%
{\mathbb{E}\left(  \left\vert \mathcal{I}_{n}^{(i)}(f)\right\vert ^{m}\right)
^{\frac{1}{m}}}\leq c(m)$, which yields the conclusion.\hfill\hbox{\vrule
height 5pt width 5pt depth 0pt}\medskip\newline

\subsection{Fluctuations of local random
fields\label{fluctuationslocalrandomfields}}

As mentioned in Section \ref{section:statementresults}, when $N$ goes to
infinity, the fields $(V^{N}_{n})_{n\geq0}$ converge in distribution to a
sequence of independent centered Gaussian random fields $(V_{n})_{n \geq0}$
whose covariances are characterized by
\[%
\begin{array}
[c]{l}%
C_{V_{n}}(f,\phi):=\mathbb{E}( V_{n}(f) V_{n}(\phi) )=\mbox{\rm Cov}_{n,\eta
_{n-1}}(f,\phi),
\end{array}
\]
for any $f, \phi\in\mathcal{B }_{b}(E)$.

We recall that for any $n \geq1$, $q \geq1$, and any $q-$tensor product
function
\[
f=\otimes_{1\leq i\leq q}f_{i} \in\mbox{\rm Osc}(E)^{\otimes q},
\]
the $q$-moments of a centered Gaussian random field $V$ are given by the Wick
formula
\begin{equation}
\label{Wick}\mathbb{E}\left(  V^{\otimes q}(f) \right)  = \sum_{\boldsymbol{i}%
\in\pi(q)} \prod_{1\leq\ell\leq q/2} \mathbb{E}(V(f_{i_{2\ell-1}})
V(f_{i_{2\ell}})),
\end{equation}
where $\pi(q)$ denotes the set of pairings of $\{ 1,\ldots, q\}$, i.e. the set
of partitions $\boldsymbol{i}$ of $\{ 1,\ldots, q\}$ into pairs
$\boldsymbol{i}_{1}=\{i_{1},i_{2}\},\ldots,\boldsymbol{i}_{q/2}=\{i_{q-1}%
,i_{q}\}$. Notice that when $q$ is odd, both sides of the above formula are
equal to zero.

In the following, we give quantitative bounds on the convergence speed for
product-form functionals of the fields $V_{n}^{N}$.

\begin{prop}
\label{borne-produit} One has the following bound, valid for any
$f=(f_{i})_{1\leq i\leq p}\in\mbox{\rm Osc}(E)^{p}$, integers $a=(a_{i}%
)_{1\leq i\leq p}$, $n\geq0$ and $N\geq1$:
\[
\left|  \mathbb{E}( V_{a_{1}}^{N}(f_{1}) \cdots V_{a_{p}}^{N}(f_{p}) ) -
\mathbb{E}( V_{a_{1}}(f_{1}) \cdots V_{a_{p}}(f_{p}) )\right|  \leq c(p)
/\sqrt{N}.
\]

\end{prop}

To prove the proposition, we use the following lemma.

\begin{lem}
\label{controle-prod} Consider a sequence of $N$ independent random variables
$(Z_{i})_{1\leq i\leq N}$ with distributions $\left(  \mu_{i}\right)  _{1\leq
i\leq N}$ on $E$, and define the empirical random fields $V^{N}$ for
$f\in\mbox{\rm Osc}(E)$ by
\[
V^{N}(f):=N^{-1/2}\sum_{j=1}^{N} (f(Z_{j}) - \mu_{j}(f)).
\]
Finally, let $\overline{V}^{N}$ denote a centered Gaussian random field with
covariance function defined for any $f, \phi\in\mbox{\rm Osc}(E)$ by
\begin{align*}
C_{\overline{V}^{N}}(f,\phi)  &  =\mathbb{E}\left(  \overline{V}^{N}(f)
\overline{V}^{N}(\phi) \right)  =\frac{1}{N}\sum_{i=1}^{N}\mbox{\rm cov}_{\mu
_{i}}(f,\phi)
\end{align*}
where
\[
\mbox{\rm cov}_{\mu_{i}}(f,\phi):=\mu_{i}\left(  \left[  f-\mu_{i}(f)\right]
\left[  \phi-\mu_{i}(\phi)\right]  \right)  .
\]
For any $1\leq q\leq N$, and any $q-$tensor product function
\[
f=\otimes_{1\leq i\leq q}f_{i} \in\mbox{\rm Osc}(E)^{\otimes q},
\]
one has that
\begin{equation}
\label{estimationCovV}\left|  \mathbb{E}\left(  \left[  V^{N}\right]
^{\otimes q}(f) \right)  - \mathbb{E}\left(  \left[  \overline{V}^{N}\right]
^{\otimes q}(f) \right)  \right|  \leq c(q) \times N^{-\rho(q)},
\end{equation}
where $\rho(q):=1$ for even $q$, and $\rho(q):=1/2$ for odd $q$.
\end{lem}

\noindent\textbf{Proof:}\newline We write
\[
V^{N}(f_{i})=\frac{1}{\sqrt{N}}\sum_{1\leq j\leq N}f^{(j)}_{i}(Z_{j}%
)\quad\mbox{\rm with}\quad f^{(j)}_{i}=f_{i}-\mu_{j}(f_{i}).
\]
Expanding the product, we get that
\[
N^{q/2}~\mathbb{E} \left(  \left[  V^{N}\right]  ^{\otimes q}(f) \right)  =
\sum_{1 \leq j_{1},\ldots, j_{q} \leq N} \mathbb{E}( f^{(j_{1})}_{1}(Z_{j_{1}%
}) \cdots f^{(j_{q})}_{q}(Z_{j_{q}}) ).
\]
Each term in the above r.h.s. such that an index $j_{i}$ appears exactly once
in the list $(j_{1},\ldots, j_{q})$ must be zero, so the only terms that may
contribute to the sum are those for which every index appears at least twice.
In the case where $q$ is odd, the number of such combinations of indices is
bounded above by $c(q) N^{(q-1)/2}$, for some finite constant $c(q)<\infty$
depending only on $q$. Since each expectation is bounded in absolute value by
1, we are done.

Now assume that $q$ is even. Consider a pairing $\boldsymbol{i}$ of
$\{1,\ldots, q \}$ given by $\boldsymbol{i}_{1}=\{i_{1},i_{2}\},\ldots
,\boldsymbol{i}_{q/2}=\{i_{q-1},i_{q}\}$, and a combination of indices
$j_{1},\ldots, j_{q}$ such that $j_{a} = j_{b}$ whenever $a,b$ belong to the
same pair, while $j_{a} \neq j_{b}$ otherwise. Denoting by $k_{r}$ the value
of $j_{a}$ when $a \in\boldsymbol{i}_{r}$, and using independence, we see that
the contribution of this combination to the sum is
\[
\mathbb{E}( f^{(j_{1})}_{1}(Z_{j_{1}}) \cdots f^{(j_{q})}_{q}(Z_{j_{q}}) ) =
\mbox{\rm cov}_{\mu_{k_{1}}}(f_{i_{1}}, f_{i_{2}}) \cdots\mbox{\rm cov}_{\mu
_{k_{q/2}}}(f_{i_{q-1}}, f_{i_{q}}).
\]
Every combination of indices in which every index appears exactly twice is of
the form we have just described. Then, the number of combinations in which
every index appears at least twice, but that are not of the previous form, is
$O(N^{q/2-1})$. As a consequence
\[%
\begin{array}
[c]{l}%
N^{q/2}~\mathbb{E}( \left(  V^{N}\right)  ^{\otimes q}(f) )\\
\\
= \displaystyle\sum_{\boldsymbol{i}\in\pi(q)}~\sum_{k \in\langle q/2,N\rangle}
\mbox{\rm cov}_{\mu_{k_{1}}}(f_{i_{1}}, f_{i_{2}}) \cdots\mbox{\rm cov}_{\mu
_{k_{q/2}}}(f_{i_{q-1}}, f_{i_{q}}) +\mbox{\rm O}\left(  N^{q/2-1}\right)  ,
\end{array}
\]
where $\langle p,N\rangle$ stands for the set of all $(N)_{p}=N!/(N-p)!$
one-to-one mappings from $[p]:=\{1,\ldots,p\}$ into $[N]$. On the other hand,
for any function $\varphi\in\mathbb{R}^{ [N]^{[p]} }$ such that $|\varphi|
\leq1$, we have
\[
\left\vert \frac{1}{(N)_{p}}\sum_{ k \in\langle p,N\rangle}\varphi(k)-
\frac{1}{N^{p}}\sum_{k \in[N]^{[p]}}\varphi(k)\right\vert \leq(p-1)/N
\]
(a detailed proof of this formula is provided in Proposition 8.6.1
in~\cite{DMo04}). Now note that
\[%
\begin{array}
[c]{l}%
\sum_{\boldsymbol{i}\in\pi(q)}~\frac{1}{N^{q/2}}\sum_{k\in[N]^{[q/2]}}
\mbox{\rm cov}_{\mu_{k_{1}}}(f_{i_{1}}, f_{i_{2}}) \cdots\mbox{\rm cov}_{\mu
_{k_{q/2}}}(f_{i_{q-1}}, f_{i_{q}} )\\
\\
= \sum_{\boldsymbol{i}\in\pi(q)}~\prod_{1\leq\ell\leq q/2} \frac{1}{N}%
\sum_{1\leq j \leq N}\mbox{\rm cov}_{\mu_{j}}(f_{i_{2\ell-1}}, f_{i_{2 \ell}%
})\\
\\
= \sum_{\boldsymbol{i}\in\pi(q)}~\prod_{1\leq\ell\leq q/2} C_{\overline{V}%
^{N}}(f_{i_{2\ell-1}}, f_{i_{2 \ell}})= \mathbb{E}\left(  \left(  \overline
{V}^{N}\right)  ^{\otimes q}(f) \right)  ,
\end{array}
\]
where the last identity uses the Wick formula \eqref{Wick}.

This yields that
\[
N^{q/2}~\mathbb{E}( \left(  V^{N}\right)  ^{\otimes q}(f) ) =\left(  N\right)
_{q/2}\mathbb{E}\left(  \left(  \overline{V}^{N}\right)  ^{\otimes q}(f)
\right)  +\mbox{\rm O}\left(  N^{q/2-1}\right)
\]
We end the proof of (\ref{estimationCovV}) using the fact that $0\leq\left(
1-(N)_{p}/N^{p}\right)  \leq(p-1)^{2}/N$, for any $p\leq N$. This ends the
proof of the lemma. \hfill\hbox{\vrule height 5pt width 5pt depth 0pt}\medskip
\newline

\begin{lem}
Given an even number $q$ and a collection of functions $(f_{i})_{1\leq i\leq
q}\in\mbox{\rm
Osc}(E)^{q}$, for any $n\geq0$ and $N\geq1$, we have
\begin{equation}%
\begin{array}
[c]{l}%
\left\Vert \prod_{1\leq\ell\leq q/2}\mbox{\rm Cov}_{n,\eta_{n-1}^{N}} (f_{2
\ell-1},f_{2 \ell})-\prod_{1\leq\ell\leq q/2}\mbox{\rm Cov}_{n,\eta_{n-1}}
(f_{2 \ell-1},f_{2 \ell}) \right\Vert _{m}\\
\\
\leq c(q,m)/\sqrt{N}.
\end{array}
\label{cov-estimate}%
\end{equation}

\end{lem}

\noindent\mbox{\bf Proof:}\newline Combining (\ref{Lm-estimate}) and
(\ref{H2}) with the generalized Minkowski inequality, we obtain that, for any
$f, \phi\in\mathrm{Osc}(E)$,
\begin{equation}
\sqrt{N}~\left\Vert \mbox{\rm Cov}_{n,\eta_{n-1}^{N}}(f,\phi
)-\mbox{\rm Cov}_{n,\eta_{n-1}}(f,\phi)\right\Vert _{m}\leq c^{^{\prime}}(m).
\label{Lm-estimate-cov}%
\end{equation}
We end the proof of (\ref{cov-estimate}) using the bound
\[
\left\vert \prod_{1\leq i\leq m}u_{i}-\prod_{1\leq i\leq m}v_{i} \right\vert
\leq\sup(|u_{i}|, |v_{i}| ; \ 1 \leq i \leq m)^{m-1} \sum_{1\leq i \leq m}
\vert u_{i} - v_{i} \vert,
\]
valid for all $u=(u_{i})_{1\leq i\leq m}\in\mathbb{R}^{m}$ and any
$v=(v_{i})_{1\leq i\leq m}\in\mathbb{R}^{m}$.

\hfill\hbox{\vrule height 5pt width 5pt depth 0pt}\medskip\newline

We now come to the proof of Proposition \ref{borne-produit}.

\textbf{Proof of Proposition \ref{borne-produit}:}\newline Assume that the
$a_{i}$ are ordered so that $a_{1}\leq\ldots\leq a_{\ell}<a_{\ell+1}%
=\cdots=a_{\ell+q}$, where $\ell+q=p$. Set
\[
A^{N}:=V_{a_{1}}^{N}(f_{1})\cdots V_{a_{\ell}}^{N}(f_{\ell})\quad
\mbox{\rm and}\quad B^{N}:=V_{a}^{N}(f_{\ell+1})\cdots V_{a}^{N}(f_{\ell+q})
\]
where $a:=a_{p}$. Given $\mathcal{F}_{a-1}^{N}$, we let $\overline{V}_{a}^{N}$
be a sequence of Gaussian random fields with covariance function defined for
any $f, \phi\in\mbox{\rm Osc}(E)$ by
\[
C_{\overline{V}_{a}^{N}}(f,\phi)=\mbox{\rm Cov}_{a,\eta_{a-1}^{N}}(f,\phi)
\]
and we set
\[
\overline{B}^{N}:=\overline{V}_{a}^{N}(f_{\ell+1})\cdots\overline{V}_{a}%
^{N}(f_{\ell+q})\quad\mbox{\rm and}\quad B:=V_{a}(f_{\ell+1})\cdots
V_{a}(f_{\ell+q})\
\]
Now $\mathbb{E}(A^{N}B^{N})=\mathbb{E}(A^{N}\times\mathbb{E}(B^{N}%
|\mathcal{F}_{a-1}^{N}))$, and, by Lemma \ref{controle-prod}, one has the
deterministic bound
\[
\left\vert \mathbb{E}(B^{N}|\mathcal{F}_{a-1}^{N})-\mathbb{E}\left(
\overline{B}^{N}\left\vert \mathcal{F}_{a-1}^{N}\right.  \right)  \right\vert
\leq c(q)/\sqrt{N}%
\]

On the other hand, combining (\ref{cov-estimate}) with Wick's formula
\eqref{Wick}
\[
\mathbb{E}\left(  \overline{B}^{N}\left\vert \mathcal{F}_{a-1}^{N}\right.
\right)  =\sum_{\boldsymbol{i}\in\pi(q)}\prod_{1\leq r \leq q/2}%
\mbox{\rm Cov}_{a,\eta_{a-1}^{N}}(f^{(\ell+2r - 1)},f^{(\ell+2 r)}),
\]
we deduce that
\[
\sqrt{N}~\left\Vert \mathbb{E}\left(  \overline{B}^{N}\left\vert
\mathcal{F}_{a-1}^{N}\right.  \right)  -\mathbb{E}\left(  B\right)
\right\Vert _{m}\leq c(m).
\]
Using the decomposition
\[
\mathbb{E}(A^{N}B^{N})-\mathbb{E}(A^{N}\mathbb{E}(B))=\mathbb{E}\left(
A^{N}\times\left[  \mathbb{E}(B^{N}|\mathcal{F}_{a-1}^{N})-\mathbb{E}%
(B)\right]  \right)
\]
we conclude that%

\[
\left\vert \mathbb{E}(A^{N}B^{N})-\mathbb{E}(A^{N})~\mathbb{E}(B)\right\vert
\leq c^{^{\prime}}(q)/\sqrt{N}.
\]
One then concludes by iterating the argument. \hfill
\hbox{\vrule height 5pt width 5pt depth 0pt}\medskip\newline

\subsection{Expansion of the particle estimates continued}

\label{puteverythingtogether}

We now plug the expansions obtained in Section \ref{decomp-de-W} into the
development obtained in (\ref{develop-log}), which leads, after some
rearrangement, to the following.

\begin{prop}
\label{prop-WN} For any $n\geq0$, $N\geq1$, we have the second order
decomposition
\begin{equation}%
\begin{array}
[c]{l}%
\displaystyle\frac{1}{\sqrt{N}}\sum_{0\leq q<n}W_{q}^{N}(\overline{G}%
_{q})-\frac{1}{2N}~\sum_{0\leq q<n}W_{q}^{N}(\overline{G}_{q})^{2}\\
\\
=\displaystyle\frac{1}{\sqrt{N}}\sum_{0\leq q<n}V_{q}^{N}(\overline{Q}%
_{q,n}(1))\\
\\
\hskip1cm\displaystyle-\frac{1}{2N}\sum_{0\leq k\leq p<n}\left[  V_{k}%
^{N}(\overline{Q}_{k,p+1}(1)-\overline{Q}_{k,p}(1))~V_{k}^{N}(\overline
{Q}_{k,p+1}(1)+\overline{Q}_{k,p}(1))~\right] \\
\hskip5cm\displaystyle-\frac{1}{N}U_{n}^{N}-\frac{1}{2N}Y_{n}^{N}+\frac
{1}{\sqrt{N}}\left(  \frac{n}{N}\right)  C_{2}(n,N)
\end{array}
\label{decomp-plug}%
\end{equation}
with the centered random variables
\[
U_{n}^{N}:=\sum_{0\leq k\not =l\leq q<p<n}V_{k}^{N}\left(  d_{k,q}%
(\overline{G}_{q})\right)  V_{l}^{N}\left(  d_{l,p}(\overline{G}_{p})\right)
\]

\[
Y_{n}^{N}:=\displaystyle\sum_{0\leq k<l\leq q<n}V_{k}^{N}\left[
d_{k,q}(\overline{G}_{q})\right]  V_{l}^{N}\left[  d_{l,q}(\overline{G}%
_{q})\right]
\]
and some remainder term such that $||C_{2}(n,N)||_{m}\leq c(m)$, for all
$m\geq1$.
\end{prop}

\noindent\mbox{\bf Proof:}\newline

By Theorem \ref{theo-key}, we may replace $W_{q}^{N}$ by $W_{q}^{N}$ in the
linear terms of the expression we want to expand, i.e. the l.h.s. of
(\ref{decomp-plug}), while committing at most an error of the form
\[
\frac{1}{\sqrt{N}}\left(  \frac{n}{N}\right)  C_{3}(n,N),
\]
where for all $m\geq1$
\[
||C_{3}(n,N)||_{m}\leq c(m).
\]

On the other hand, using the cruder expansion provided by Theorem
\ref{decomp-exacte}, we may replace $W_{q}^{N}$ by just $\sum_{p=0}^{q}%
V_{p}^{N}\left[  d_{p,q}(\overline{G}_{q})\right]  $ in the quadratic terms
appearing in the l.h.s. of (\ref{decomp-plug}), and commit an overall error of
the form
\[
\frac{1}{\sqrt{N}}\left(  \frac{n}{N}\right)  C_{4}(n,N),
\]
where for all $m\geq1$
\[
||C_{4}(n,N)||_{m}\leq c^{^{\prime}}(m).
\]

By the definition of $W_{q}^{N}$ given in (\ref{2dec}), we have
\begin{align*}
W_{q}^{N}(\overline{G}_{q})  &  =\displaystyle\sum_{p=0}^{q}V_{p}^{N}\left[
d_{p,q}(\overline{G}_{q})\right] \\
&  -\displaystyle\frac{1}{\sqrt{N}}\sum_{0\leq p<q}~\left[  \sum_{k=0}%
^{p}V_{k}^{N}\left[  d_{k,p}(\overline{G}_{k})\right]  \right]  \left[
\sum_{k=0}^{p}V_{k}^{N}\left[  d_{k,q}(\overline{G}_{q})\right]  \right]
\end{align*}
so that
\[%
\begin{array}
[c]{l}%
\displaystyle\frac{1}{\sqrt{N}}\sum_{0\leq q<n}W_{q}^{N}(\overline{G}_{q})\\
\\
=\displaystyle\frac{1}{\sqrt{N}}\sum_{0\leq p<n}V_{p}^{N}\left[  \sum_{p\leq
q<n}d_{p,q}(\overline{G}_{q})\right] \\
\\
\hskip2cm-\displaystyle\frac{1}{N}\sum_{0\leq q<n}\sum_{0\leq p<q}~\left[
\sum_{k=0}^{p}V_{k}^{N}\left[  d_{k,p}(\overline{G}_{k})\right]  \right]
\left[  \sum_{k=0}^{p}V_{k}^{N}\left[  d_{k,q}(\overline{G}_{q})\right]
\right]  .
\end{array}
\]
We recall that
\[
\sum_{p\leq q<n}d_{p,q}(\overline{G}_{q})=\sum_{p\leq q<n}\left[  \overline
{Q}_{p,q+1}(1)-\overline{Q}_{p,q}(1)\right]  =\overline{Q}_{p,n}(1)-1
\]
so that on the one hand we have
\[
\sum_{0\leq p<n}V_{p}^{N}\left[  \sum_{p\leq q<n}d_{p,q}(\overline{G}%
_{q})\right]  =\sum_{0\leq p<n}V_{p}^{N}\left[  \overline{Q}_{p,n}(1)\right]
,
\]
whereas, on the other hand, we have
\[%
\begin{array}
[c]{l}%
\sum_{0\leq p<q<n}\displaystyle\left[  \sum_{k=0}^{p}V_{k}^{N}\left[
d_{k,p}(\overline{G}_{k})\right]  \right]  \left[  \sum_{k=0}^{p}V_{k}%
^{N}\left[  d_{k,q}(\overline{G}_{q})\right]  \right] \\
\\
=\displaystyle\sum_{0\leq k\leq p<q<n}V_{k}^{N}\left[  d_{k,p}(\overline
{G}_{k})\right]  V_{k}^{N}\left[  d_{k,q}(\overline{G}_{q})\right]  +U_{n}%
^{N}\\
\\
=\displaystyle\sum_{0\leq k<q<n}V_{k}^{N}\left[  \sum_{k\leq p<q}%
d_{k,p}(\overline{G}_{k})\right]  V_{k}^{N}\left[  d_{k,q}(\overline{G}%
_{q})\right]  +U_{n}^{N}\\
\\
=\displaystyle\sum_{0\leq k<q<n}V_{k}^{N}\left[  \overline{Q}_{k,q}(1)\right]
V_{k}^{N}\left[  d_{k,q}(\overline{G}_{q})\right]  +U_{n}^{N}.
\end{array}
\]
This implies that
\[%
\begin{array}
[c]{l}%
\displaystyle\frac{1}{\sqrt{N}}\sum_{0\leq q<n}W_{q}^{N}(\overline{G}_{q})\\
\\
=\displaystyle\frac{1}{\sqrt{N}}\sum_{0\leq p<n}V_{p}^{N}\left[  \overline
{Q}_{p,n}(1)\right] \\
\\
\hskip2cm\displaystyle-\frac{1}{N}\displaystyle\sum_{0\leq q<n}\sum_{0\leq
p<q}V_{p}^{N}\left[  \overline{Q}_{p,q}(1)\right]  V_{p}^{N}\left[
d_{p,q}(\overline{G}_{q})\right]  -\frac{1}{N}~U_{n}^{N}.
\end{array}
\]
It remains to analyze the quadratic part, which we write as
\[%
\begin{array}
[c]{l}%
\displaystyle\sum_{0\leq q<n}\left(  \sum_{0\leq p\leq q}V_{p}^{N}\left[
d_{p,q}(\overline{G}_{q})\right]  \right)  ^{2}=\sum_{0\leq q<n}\sum_{0\leq
p\leq q}V_{p}^{N}\left[  d_{p,q}(\overline{G}_{q})\right]  ^{2}+Y_{n}^{N}.
\end{array}
\]
Now notice that
\[%
\begin{array}
[c]{l}%
-\displaystyle\sum_{0\leq p<q}V_{p}^{N}\left[  \overline{Q}_{p,q}(1)\right]
V_{p}^{N}\left[  d_{p,q}(\overline{G}_{q})\right]  -\frac{1}{2}%
\displaystyle\sum_{0\leq p\leq q}V_{p}^{N}\left[  d_{p,q}(\overline{G}%
_{q})\right]  ^{2}\\
\\
=\displaystyle-\frac{1}{2}V_{q}^{N}\left[  d_{q,q}(\overline{G}_{q})\right]
^{2}-\sum_{0\leq p<q}V_{p}^{N}\left[  d_{p,q}(\overline{G}_{q})\right]
~V_{p}^{N}\left[  \frac{1}{2}~d_{p,q}(\overline{G}_{q})+\overline{Q}%
_{p,q}(1)\right] \\
\\
=\displaystyle-\frac{1}{2}V_{q}^{N}\left[  d_{q,q}(\overline{G}_{q})\right]
^{2}\\
\\
\hskip1cm-\sum_{0\leq p<q}V_{p}^{N}\left[  d_{p,q}(\overline{G}_{q})\right]
~V_{p}^{N}\left[  \frac{1}{2}\left[  \overline{Q}_{p,q+1}(1)-\overline
{Q}_{p,q}(1)\right]  +\overline{Q}_{p,q}(1)\right] \\
\\
=\displaystyle-\frac{1}{2}V_{q}^{N}\left[  \overline{Q}_{q,q+1}(1)-\overline
{Q}_{q,q}(1)\right]  ^{2}\\
\\
\hskip1cm-\frac{1}{2}\sum_{0\leq p<q}V_{p}^{N}\left[  \overline{Q}%
_{p,q+1}(1)-\overline{Q}_{p,q}(1)\right]  ~V_{p}^{N}\left[  \overline
{Q}_{p,q+1}(1)+\overline{Q}_{p,q}(1)\right]
\end{array}
\]
Recalling that $\overline{Q}_{q,q}(1)=1$, we conclude that
\[%
\begin{array}
[c]{l}%
-\displaystyle\sum_{0\leq p<q}V_{p}^{N}\left[  \overline{Q}_{p,q}(1)\right]
V_{p}^{N}\left[  d_{p,q}(\overline{G}_{q})\right]  -\frac{1}{2}%
\displaystyle\sum_{0\leq p\leq q}V_{p}^{N}\left[  d_{p,q}(\overline{G}%
_{q})\right]  ^{2}\\
\\
=-\frac{1}{2}\sum_{0\leq k\leq q}V_{k}^{N}\left[  \overline{Q}_{k,q+1}%
(1)-\overline{Q}_{k,q}(1)\right]  ~V_{k}^{N}\left[  \overline{Q}%
_{k,q+1}(1)+\overline{Q}_{k,q}(1)\right]
\end{array}
\]
\hfill\hbox{\vrule height 5pt width 5pt depth 0pt}\medskip\newline

The next step is to show that both centered terms $U_{n}^{N}$ and $Y_{n}^{N}$
yield negligible contributions in (\ref{decomp-plug}).

\begin{prop}
\label{neglige-U} For any $n\geq0$, and any $N\geq1$, we have that
\[
\mathbb{E}((U_{n}^{N})^{2})\leq c \ \left(  n+\frac{n^{2}}{\sqrt{N}}\right)
.
\]

\end{prop}

\noindent\mbox{\bf Proof:}\newline

We can write
\[%
\begin{array}
[c]{l}%
\mathbb{E}((U_{n}^{N})^{2})\\
\\
=\sum\mathbb{E}\left(  V_{k}^{N}\left(  d_{k,q}(\overline{G}_{q})\right)
V_{l}^{N}\left(  d_{l,p}(\overline{G}_{p})\right)  V_{k^{\prime}}^{N}\left(
d_{k^{\prime},q^{\prime}}(\overline{G}_{q^{\prime}})\right)  V_{l^{\prime}%
}^{N}\left(  d_{l^{\prime},p^{\prime}}(\overline{G}_{p^{\prime}})\right)
\right)
\end{array}
\]
with
\[
\sum=\sum_{0\leq k\not =l\leq q<p<n}\sum_{0\leq k^{\prime}\not =l^{\prime}\leq
q^{\prime}<p^{\prime}<n}.
\]

First consider replacing each $V_{k}^{N}$ by the corresponding $V_{k}$ in the
above expectations. By Proposition \ref{borne-produit} together with
(\ref{control-gpn}), the overall error is bounded by
\[
\frac{c(p)}{\sqrt{N}}\sum\exp(-\lambda(q-k+p-l+q^{\prime}-k^{\prime
}+p^{\prime}-l^{\prime}))\leq c^{\prime}\left(  p\right)  \frac{n^{2}}{\sqrt{N}}.
\]
Now consider the corresponding sum
\[
\sum\mathbb{E}\left(  V_{k}\left(  d_{k,q}(\overline{G}_{q})\right)
V_{l}\left(  d_{l,p}(\overline{G}_{p})\right)  V_{k^{\prime}}\left(
d_{k^{\prime},q^{\prime}}(\overline{G}_{q^{\prime}})\right)  V_{l^{\prime}%
}\left(  d_{l^{\prime},p^{\prime}}(\overline{G}_{p^{\prime}})\right)  \right)
.
\]
The only possibility to have a non-zero term is when either $k=k^{\prime}$ and
$l=l^{\prime}$ or $k=l^{\prime}$ and $k^{\prime}=l$. Restricting summation to
this subset of indices, we obtain that
\[
\sum\exp(-\lambda(q-k+p-l+q^{\prime}-k^{\prime}+p^{\prime}-l^{\prime}))\leq
c^{^{\prime}} \times n.
\]
\hfill\hbox{\vrule height 5pt width 5pt depth 0pt}\medskip\newline

With a similar argument, we also obtain the following result.

\begin{prop}
\label{neglige-Y} For any $n\geq0$, $N\geq1$, we have
\[
\mathbb{E}((Y_{n}^{N})^{2})\leq c \ \left(  n+\frac{n^{2}}{\sqrt{N}}\right)
.
\]

\end{prop}

Now, we consider the remaining term in (\ref{decomp-plug}), i.e.
\[
H_{n}^{N} := \sum_{0\leq k\leq p<n} \left(  V^{N}_{k}\left[  \overline
{Q}_{k,p+1}(1)-\overline{Q}_{k,p}(1)\right]  ~V^{N}_{k}\left[  \overline
{Q}_{k,p+1}(1)+\overline{Q}_{k,p}(1)\right]  ~\right)  ,
\]
and show that it can be replaced by its expectation up to a negligible random term.

\begin{prop}
\label{biais-determin} For any $n\geq0$, $N\geq1$, we have the following bound
:
\[
\mathbb{V}(H_{n}^{N})\leq c \ n.
\]

\end{prop}

\noindent\mbox{\bf Proof:}\newline If we set
\[
J_{k,p}:=\overline{Q}_{k,p+1}(1)-\overline{Q}_{k,p}(1)\quad\mbox{\rm and}\quad
K_{k,p}:=\overline{Q}_{k,p+1}(1)+\overline{Q}_{k,p}(1)
\]
then we find that
\[
\mathbb{E}\left(  H_{n}^{N}\right)  =\sum_{0\leq k\leq p<n}\mathbb{E}\left(
V_{k}^{N}\left[  J_{k,p}\right]  ~V_{k}^{N}\left[  K_{k,p}\right]  ~\right)
\]
whence
\[%
\begin{array}
[c]{l}%
(\mathbb{E}\left(  H_{n}^{N}\right)  )^{2}\\
\\
=\displaystyle\sum_{0\leq k\leq p<n}\sum_{0\leq k^{\prime}\leq p^{\prime}%
<n}\mathbb{E}\left(  V_{k}^{N}\left[  J_{k,p}\right]  ~V_{k}^{N}\left[
K_{k,p}\right]  ~\right)  \mathbb{E}\left(  V_{k^{\prime}}^{N}\left[
J_{k^{\prime},p^{\prime}}\right]  ~V_{k^{\prime}}^{N}\left[  K_{k^{\prime
},p^{\prime}}\right]  ~\right)
\end{array}
\]
while
\[%
\begin{array}
[c]{l}%
\mathbb{E}\left(  \left(  H_{n}^{N}\right)  ^{2}\right) \\
\\
=\displaystyle\sum_{0\leq k\leq p<n}\sum_{0\leq k^{\prime}\leq p^{\prime}%
<n}\mathbb{E}\left(  V_{k}^{N}\left[  J_{k,p}\right]  ~V_{k}^{N}\left[
K_{k,p}\right]  V_{k^{\prime}}^{N}\left[  J_{k^{\prime},p^{\prime}}\right]
~V_{k^{\prime}}^{N}\left[  K_{k^{\prime},p^{\prime}}\right]  ~\right)  .
\end{array}
\]

Observe that, whenever $k\neq k^{\prime}$, the terms in the above two sums
coincide. Therefore, it remains to bound the contribution in both sums of the
terms that have $k=k^{\prime}$. In both expressions, the corresponding sum is
bounded above in absolute value by
\[
\sum_{0\leq k\leq p,p^{\prime}<n}c^{^{\prime}} \times e^{-\lambda(p^{\prime
}-k+p-k)}\leq c^{^{\prime\prime}} \times n.
\]
This ends the proof of the proposition. \hfill\hbox{\vrule height 5pt width
5pt depth 0pt}\medskip\newline

\begin{prop}
\label{limite-biais} For any $n\geq0$, $N\geq1$, we have
\[
\mathbb{E}(H_{n}^{N})=v_{n}+\epsilon_{n}^{N}\quad\mbox{with}\quad|\epsilon
_{n}^{N}|\leq c\times {n}/{\sqrt{N}} .
\]

\end{prop}

\noindent\mbox{\bf Proof:}\newline Recalling that $\overline{Q}_{p,n}%
(1)-1=\sum_{p\leq k<n}\left(  \overline{Q}_{p,k+1}-\overline{Q}_{p,k}\right)
$, we prove that
\[%
\begin{array}
[c]{l}%
V_{p}(\overline{Q}_{p,n}(1))^{2}\\
\\
=\left(  \sum_{p\leq k<n}V_{p}\left(  \overline{Q}_{p,k+1}-\overline{Q}%
_{p,k}\right)  \right)  ^{2}\\
\\
=\sum_{p\leq k<n}V_{p}\left(  \overline{Q}_{p,k+1}-\overline{Q}_{p,k}\right)
^{2}\\
\\
\hskip1cm+2\sum_{p\leq l<n} V_{p}\left(  \sum_{p\leq k<l}\left[  \overline
{Q}_{p,k+1}-\overline{Q}_{p,k}\right]  \right)  V_{p}\left(  \overline
{Q}_{p,l+1}-\overline{Q}_{p,l}\right) \\
\\
=\sum_{p\leq l<n}V_{p}\left(  \overline{Q}_{p,l+1}-\overline{Q}_{p,l}\right)
^{2}\\
\\
\hskip4cm+2\sum_{p\leq l<n} V_{p}\left(  \overline{Q}_{p,l}\right)
V_{p}\left(  \overline{Q}_{p,l+1}-\overline{Q}_{p,l}\right)
\end{array}
\]
This yields the formula
\[
V_{p}(\overline{Q}_{p,n}(1))^{2}=\sum_{p\leq l<n} V_{p}\left(  \overline
{Q}_{p,l+1}-\overline{Q}_{p,l}\right)  ~V_{p}\left(  \overline{Q}%
_{p,l+1}+\overline{Q}_{p,l}\right)
\]
Replacing each $V_{k}^{N}$ by $V_{k}$ in the expectation of $H_{n}^{N}$, we
obtain
\[%
\begin{array}
[c]{l}%
\sum_{0 \leq p \leq l < n } \mathbb{E} \left(  V_{p} \left[  \overline
{Q}_{p,l+1}(1)-\overline{Q}_{p,l}(1) \right]  V_{p}\left[  \overline
{Q}_{p,l+1}(1)+\overline{Q}_{p,l}(1) \right]  \right) \\
\\
= \sum_{0 \leq p < n } \mathbb{E}\left(  V_{p}(\overline{Q}_{p,n}%
(1))^{2}\right)  = v_{n}.
\end{array}
\]

To control the error introduced by the replacement, we use Proposition
\ref{borne-produit}, (\ref{Qpn-diff}) and (\ref{control-gpn}), so that the
overall error can be bounded above by
\[
c \ \sum_{0\leq k\leq p<n}\frac{e^{-\lambda(p-k)}}{\sqrt{N}}\leq c^{^{\prime}}
\ \frac{n}{\sqrt{N}}%
\]
This ends the proof of the proposition. \hfill\hbox{\vrule height 5pt width
5pt depth 0pt}\medskip\newline

\subsection{Central limit theorem}

\label{centrallimittheoremproof}

This section established the proof of theorem~\ref{lognormal}. Using
Proposition (\ref{develop-log}), the decomposition (\ref{decomp-plug}), and
Propositions (\ref{neglige-U}), (\ref{neglige-Y}), (\ref{biais-determin}) and
(\ref{limite-biais}), we obtain
\[
\log\overline{\gamma}_{n}^{N}(1)=\frac{1}{\sqrt{N}}\sum_{0\leq q<n}V_{q}%
^{N}(\overline{Q}_{q,n}(1))-\frac{1}{2N}~v_{n}+\varepsilon_{n}^{N},
\]
with $\varepsilon_{n}^{N}$ going to zero in probability as $n$ goes to
infinity. Thus, to prove the theorem, it remains to show that
\[
\frac{1}{\sqrt{v_{n}}}\sum_{0\leq q<n}V_{q}^{N}(\overline{Q}_{q,n}(1))
\]
converges in distribution to a standard normal. We do so using the central
limit theorem for martingale difference arrays (see e.g. \cite{Dur,Shi}). The
martingale property just comes from the fact that, for any $q\geq0$ and any
bounded function $f_{q}$, one has
\[
\mathbb{E}\left(  V_{q}^{N}(f_{q})|\mathcal{F}_{q-1}^{N}\right)  =0\ a.s.
\]

We now have to show that
\[
\frac{1}{v_{n}}\sum_{0\leq q<n}\mathbb{E}\left(  \left[  V_{q}^{N}%
(\overline{Q}_{q,n}(1))\right]  ^{2}|\mathcal{F}_{q-1}^{N}\right)
\]
converges to $1$ in probability. One easily checks from the definition that
\[
\mathbb{E}\left(  \left[  V_{q}^{N}(\overline{Q}_{q,n}(1))\right]
^{2}|\mathcal{F}_{q-1}^{N}\right)  =\mbox{\rm Cov}_{q,\eta_{q-1}^{N}%
}(\overline{Q}_{q,n}(1),\overline{Q}_{q,n}(1))
\]
We observe that
\[
v_{n}=\sum_{0\leq q<n}\mbox{\rm Cov}_{q,\eta_{q-1}}(\overline{Q}%
_{q,n}(1),\overline{Q}_{q,n}(1))
\]
and
\[%
\begin{array}
[c]{l}%
d_{n}^{N}:=\left\vert \frac{1}{v_{n}}\sum_{0\leq q<n}\mathbb{E}\left(  \left[
V_{q}^{N}(\overline{Q}_{q,n}(1))\right]  ^{2}|\mathcal{F}_{q-1}^{N}\right)
-1\right\vert \\
\\
\leq\frac{1}{v_{n}}\sum_{0\leq q<n}\left\vert \mbox{\rm Cov}_{q,\eta_{q-1}%
^{N}}(\overline{Q}_{q,n}(1),\overline{Q}_{q,n}(1))-\mbox{\rm Cov}_{q,\eta
_{q-1}}(\overline{Q}_{q,n}(1),\overline{Q}_{q,n}(1))\right\vert
\end{array}
\]
Using (\ref{Lm-estimate-cov}), we see that
\[
\mathbb{E}(d_{n}^{N})\leq c \ \left(  \frac{n}{v_{n}}\right)  \frac{1}%
{\sqrt{N}},
\]
so we can conclude using (\ref{assumption-borneinf}).

The last point to be checked is the asymptotic negligibility condition, that
is, for all $\epsilon>0$, we have to prove that
\[
\frac{1}{v_{n}}\sum_{0\leq q<n}\mathbb{E}\left(  \left[  V_{q}^{N}%
(\overline{Q}_{q,n}(1))\right]  ^{2}{\mathord{1\!{\rm l}}}\left(  \left[
V_{q}^{N}(\overline{Q}_{q,n}(1))\right]  ^{2}\geq\epsilon\text{ }v_{n}\right)
|\mathcal{F}_{q-1}^{N}\right)
\]
goes to zero in probability. By Schwarz's inequality and (\ref{Lm-estimate}),
the expectation of this expression is bounded above by
\[
c^{^{\prime}} \ \left(  \frac{n}{v_{n}}\right)  \frac{1}{(\epsilon
v_{n})^{1/2}},
\]
This ends the proof of the theorem. \hfill\hbox{\vrule height 5pt width 5pt
depth 0pt}\medskip\newline

\end{document}